\documentclass[12pt,halfline,a4paper]{ouparticle}

\usepackage{float}
\usepackage{amsthm}
\usepackage{xcolor}
\usepackage{dsfont}
\usepackage{graphicx}
\usepackage{subcaption}
\usepackage{hyperref}

\newtheorem{theorem}{Theorem}
\newtheorem{proposition}{Proposition}

\newtheorem{ex}{Example}
\newtheorem{cor}{Corollary}

\usepackage{mathtools}

\DeclarePairedDelimiter\floor{\lfloor}{\rfloor}

\begin{document}

\title{Picturesque convolution-like recurrences and partial sums' generation}

\author{%
\name{Ignas Gasparavičius, Andrius Grigutis, Juozas Petkelis}
\address{Institute of Mathematics, Vilnius University\\
Naugarduko g. 24, LT-03225, Vilnius, Lithuania}
\email{ignas.gasparavicius@mif.stud.vu.lt, andrius.grigutis@mif.vu.lt, juozas.petkelis@mif.stud.vu.lt}
}

\abstract{Let ${\pmb b}=\{b_0,\,b_1,\,\ldots\}$ be the known sequence of numbers such that $b_0\neq0$. In this work, we develop methods to find another sequence ${\pmb a}=\{a_0,\,a_1,\,\ldots\}$ that is related to ${\pmb b}$ as follows: $a_n=a_0\,b_{n+m}+a_1\,b_{n+m-1}+\ldots+a_{n+m}\,b_0$, $n\in\mathbb{N}\cup\{0\}$, $m\in\mathbb{N}$. We show the connection of $\lim_{n\to\infty}a_n$ with $a_0,\,a_1,\,\ldots,\,a_{m-1}$ and provide varied examples of finding the sequence ${\pmb a}$ when ${\pmb b}$ is given. We demonstrate that the sequences ${\pmb a}$ may exhibit pretty patterns in the plane or space. Also, we show that the properly chosen sequence ${\pmb b}$ may define {\pmb a} as some famous sequences, such as the partial sums of the Riemann zeta function, etc.
}

\date{\today}

\keywords{
linear recurrence, initial values, Maclaurin series, power series, partial sums, Riemann hypothesis
}

\maketitle

\section{Introduction}\label{sec1}

\subsection{Short overview of the recurrence phenomena}

Recurrence relations play an important role in the surrounding phenomena. For instance, in computer science, they serve as fundamental tools for modeling iterative and recursive processes where the current state depends on previous states. These mathematical structures arise in various applications, from algorithm design and data structures to artificial intelligence and signal processing.  One notable class of recurrence-based methods involves convolution-like recurrent sequences. For example, recurrent sequences filter signals in signal processing, removing noise and enhancing essential features. These filters use a recurrence relation where the current output depends not only on the current and past input values but also on past output values. In machine learning, convolutional neural networks (CNNs) rely on discrete convolutions, and recurrent neural network (RNN) models utilize recurrence relations to predict future values based on historical data.
These sequences provide a robust framework for understanding how data evolves over time and how patterns can be extracted from seemingly chaotic inputs. By leveraging the properties of recurrence relations, people develop efficient algorithms and systems capable of handling various types of data. For a broader description of the mentioned facts and more, see, for instance, \cite{greene}, \cite{enchantments} \cite{Tucker}, \cite{penguin}, \cite{Wilf}.

\subsection{Introduction to the article}

In this work, we suppose that the sequence of complex numbers ${\pmb b}:=\{b_0,\,b_1,\,\ldots\}$ is known. In addition, we assume that $b_0\neq 0$. We fix $m\in\mathbb{N}$ and study the convolution-like recurrence 
\begin{align}\label{eq:seq}
a_n=\sum_{j=0}^{n+m}b_{n+m-j}\,a_j,\,n\in\mathbb{N}_0:=\mathbb{N}\cup\{0\}.
\end{align}

We provide an algorithm to compute ${\pmb a} = \{a_0,\,a_1,\,\ldots\}$. As we will soon see, the mentioned algorithm is based on the auxiliary sequences that connect $a_n$ with the linear combination of $a_0,\,a_1,\,\ldots,\,a_{m-1}$ for all $n\in\mathbb{N}_0$. More precisely, due to the assumption $b_0\neq 0$, the recurrence \eqref{eq:seq} implies
\begin{align}\label{eq:seq_v1}
a_n=\frac{1}{b_0}\left(a_{n-m}-\sum_{j=0}^{n-1}b_{n-j}a_j\right),\,n=m,\,m+1,\,\ldots,\,m\in\mathbb{N}.
\end{align}

Similarly as in \eqref{eq:seq_v1}, we define $m\in\mathbb{N}$ recurrent sequences $\alpha_0(n)$, $\alpha_1(n)$, $\ldots$, $\alpha_{m-1}(n)$:
\begin{align}\label{rec_for_expres}
\begin{cases}
\alpha_0(n)=\frac{1}{b_0}\left(\alpha_0(n-m)-\sum\limits_{j=0}^{n-1}b_{n-j}\alpha_{0}(j)\right)\\
\alpha_1(n)=\frac{1}{b_0}\left(\alpha_1(n-m)-\sum\limits_{j=0}^{n-1}b_{n-j}\alpha_{1}(j)\right)\\
\,\,\vdots\\
\alpha_{m-1}(n)=\frac{1}{b_0}\left(\alpha_{m-1}(n-m)-\sum\limits_{j=0}^{n-1}b_{n-j}\alpha_{m-1}(j)\right)
\end{cases},\,
n=m+1,\,m+2,\,\ldots,
\end{align}
whose initial values are given in Table \ref{T}. 
\begin{table}[H]
\centering
\begin{tabular}{|c||c|c|c|c|c|c|} 
 \hline
$n$ & $\alpha_0(n)$ & $\alpha_1(n)$ & $\alpha_2(n)$ & $\ldots$& $\alpha_{m-2}(n)$&$\alpha_{m-1}(n)$\\ 
 \hline\hline
 $0$ & $1$ & $0$ & $0$ & $\ldots$&$0$&$0$\\ 
 \hline
 $1$ & $0$ & $1$ & $0$ & $\ldots$&$0$&$0$\\
 \hline
 $2$ & $0$ & $0$ & $1$ & $\ldots$&$0$&$0$\\
 \hline
 $\vdots$ & $\vdots$ & $\vdots$ & $\vdots$&$\ddots$&$\vdots$&$\vdots$ \\
 \hline
 $m-1$ & $0$ & $0$ & $0$ &$\ldots$&$0$&$1$\\ 
 \hline
 $m$ & $\frac{1-b_m}{b_{0}}$ & $-\frac{b_{m-1}}{b_0}$ & $-\frac{b_{m-2}}{b_0}$ &$\ldots$&$-\frac{b_2}{b_0}$&$-\frac{b_1}{b_0}$\\ 
 \hline
\end{tabular}
\caption{The initial values for the recurrent sequences $\alpha_0(n),\,\alpha_1(n),\,\ldots,\,\alpha_{m-1}(n)$ in \eqref{rec_for_expres}.}
\label{T}
\end{table}
{\sc Note 1:} {\it We emphasize that the penultimate row and column in Table \ref{T} are only present if $m\in\{2,\,3,\,\ldots\}$.}

The following Proposition provides the mentioned expressiveness of $a_n,\,n\in\mathbb{N}_0$ via $a_0,\,a_1,\,\ldots,\,a_{m-1}$.

\begin{proposition}\label{lem}
Let $m\in\mathbb{N}$. For the defined sequences $\alpha_0(n),\,\alpha_1(n),\,\ldots,\,\alpha_{m-1}(n)$ it holds that
\begin{align}\label{a_n_expr}
a_n=\alpha_0(n)a_0+\alpha_1(n)a_1+\ldots+\alpha_{m-1}(n)a_{m-1},\,n\in\mathbb{N}_0.
\end{align}
\end{proposition}

\begin{proof}
See Section \ref{sec:proofs}.
\end{proof}

Eq. \eqref{a_n_expr} suggests that the initial values $a_0$, $a_1$, $\ldots$, $a_{m-1}$, $m\in\mathbb{N}$ for the recurrence \eqref{eq:seq_v1} can be expressed via $\lim_{n\to\infty}a_n$ if the limit behavior, as $n\to\infty$, of the right hand-side of \eqref{a_n_expr} is known. In the next section (Theorem \ref{thm}), we prove that the values of sequences $\alpha_0(n),\,\alpha_1(n),\,\ldots,\,\alpha_{m-1}(n)$ can be computed as coefficients of certain power series. We provide conditions (Theorem \ref{thm:sum_b_to_1}) when every limit of $\alpha_0(n),\,\alpha_1(n),\,\ldots,\,\alpha_{m-1}(n)$, as $n\to\infty$, exists and provide the expressions of these limits. In Theorem \ref{thm:expresions_general} and Corollary \ref{cor:expressions}, we show that every single term $a_0$, $a_1$, $\ldots$, $a_{m-1}$ can be expressed via $\lim_{n\to\infty}a_n$. As an important application of this study, in Propositions \ref{ex_zeta_1}--\ref{proposition_pi}, we show that the sequence ${\pmb b}$ can be chosen in a way that its outputs \eqref{rec_for_expres} generate partial sums that approximate some famous constants (values of the Riemann zeta-function, $\pi$, $e$) and these constants can be computed as $b_1+2b_2+3b_3+\ldots$ It is intriguing that sometimes these $\alpha$-sequences \eqref{rec_for_expres} (${\pmb a}$ as well) may exhibit striking patterns in space or plane, as seen in Figures \ref{F1}--\ref{F7} below. In the last Section \ref{sec:examples}, we provide many and varied examples demonstrating what $\alpha$-sequences \eqref{rec_for_expres} are implied when a certain sequence ${\pmb b}$ is chosen. In short, the essence of this work revolves around finding the Maclaurin series (the Taylor series at zero) of rational or irrational functions and utilizing some ingenuity with their coefficients.

\subsection{Prior work}

To our knowledge, the expressions of type \eqref{rec_for_expres} with $m=2$ and ${\pmb b}$ representing the pro\-ba\-bi\-li\-ty distribution were first introduced in \cite{JJ}. The authors there constructed the pro\-ba\-bi\-li\-ty distribution function $\mathbb{P}(\max\{0,\,X_1,\,X_1+X_2,\,\ldots\}<x),\,x>0$, where $X_1,\,X_2,\,\ldots$ are certain random variables. In \cite{AJ}, it was proven the connection between recurrences of type \eqref{rec_for_expres} (with $m=2$ and ${\pmb b}$ representing the pro\-ba\-bi\-li\-ty distribution) and the probability generating functions. Then, the setup of the distribution function $\mathbb{P}(\max\{X_1,\,X_1+X_2,\,\ldots\}<x)$ based on generating functions was gradually developed by the second-named author of this work and coauthors, see \cite{A2}, \cite{A1}, \cite{AA} and references therein. In the mentioned works, the sequence ${\pmb b}$ was always considered as a discrete probability distribution, and the sequence ${\pmb a}$, which was desired to find, was the distribution function of $\max\{X_1,\,X_1+X_2,\,\ldots\}$. In the described probabilistic context, the probability mass function and the distribution function of $\max\{X_1,\,X_1+X_2,\,\ldots\}$ are related similarly as \eqref{eq:seq}. In this work, by letting ${\pmb b}$ be an arbitrary sequence, we depart from the comfortable probabilistic context and examine the generalizations of the previously used methods.

\section{Main results}

In this section, we formulate the main results of our work, demonstrate some of its applications, and depict several graphs of the considered sequences.

Let $s\in\mathbb{C}$ and define the following formal power series for the previously introduced sequences ${\pmb a}$, ${\pmb b}$,  and $\alpha_0(n),\,\alpha_1(n),\,\ldots,\,\alpha_{m-1}(n)$, $n\in\mathbb{N}_0$, $m\in\mathbb{N}$:
\begin{align*}
&A(s):=\sum_{n=0}^{\infty}a_n s^n,\,B(s):=\sum_{n=0}^{\infty}b_n s^n,\\
&G\alpha_0(s):=\sum_{n=0}^{\infty}\alpha_0(n) s^n,\,
G\alpha_1(s):=\sum_{n=0}^{\infty}\alpha_1(n) s^n,\,\ldots,\,
G\alpha_{m-1}(s):=\sum_{n=0}^{\infty}\alpha_{m-1}(n) s^n.
\end{align*}

For the defined power series and considered sequences, the following statement is correct.

\newpage

\begin{theorem}\label{thm}
Let ${\pmb b}$ be the known sequence of numbers such that $b_0\neq0$. If the sequence ${\pmb a}$ is related to ${\pmb b}$ as \eqref{eq:seq}, then, for $m\in\mathbb{N}$ and some $s\in\mathbb{C}$, following relations hold:
\begin{align}
\sum_{k=0}^{m-1}a_k\sum_{n=k}^{m-1}b_{n-k}\,s^n&=A(s)(B(s)-s^m),\label{a_over_b}\\
A(s)&=\sum_{j=0}^{m-1}a_j\,G\alpha_j(s),\label{a_over_alpha}\\
G\alpha_k(s)(B(s)-s^m)&=\sum_{n=k}^{m-1}b_{n-k}\,s^n, \text{ for every } k=0,\,1,\,\ldots,\,m-1,\label{alpha_over_b_with_initial}\\
a_n-a_{n-1}&=\sum_{k=0}^{m-1}a_k(\alpha_k(n)-\alpha_k(n-1)),\,n\in\mathbb{N},\label{a_diff}\\
\alpha_k(n)-\alpha_k(n-1)&=
\frac{1}{n!}
\lim_{s\to0}
\left(
\frac{\partial^n}{\partial s^n}
(1-s)G\alpha_k(s)
\right)
,\,n\in\mathbb{N},\label{alpha_diff}\\
\alpha_k(n)&=\frac{1}{n!}
\lim_{s\to0}
\left(
\frac{\partial^n}{\partial s^n}
G\alpha_k(s)
\right)
,\,n\in\mathbb{N}_0.\label{alpha_ne_diff}
\end{align}
\end{theorem}

\begin{proof}
See Section \ref{sec:proofs}.
\end{proof}

Eqs. \eqref{alpha_diff} and \eqref{alpha_ne_diff} of Theorem \ref{thm} imply the following corollary.

\begin{cor}\label{Tuite}
Under the assumptions of Theorem \ref{thm}, it holds that
\begin{align}
\limsup_{n\to\infty}\left|\alpha_k(n)-\alpha_k(n-1)\right|^{1/n}=\frac{1}{R_k} \text{ for every } k=0,\,1,\,\ldots,\,m-1,\label{R}
\end{align}
where $R_k>0$ correspondingly are the supremums of such positive numbers $r_k$ that every $(1-s)G\alpha_k(s)$ is analytic in $|s|<r_k$. Moreover,
\begin{align}
\limsup_{n\to\infty}\left|\alpha_k(n)\right|^{1/n}=\frac{1}{\tilde{R}_k} \text{ for every } k=0,\,1,\,\ldots,\,m-1,\label{R_1}
\end{align}
where $\tilde{R}_k>0$ correspondingly are the supremums of such positive numbers $\tilde{r}_k$ that every $G\alpha_k(s)$ is analytic in $|s|<\tilde{r}_k$.
\end{cor}

\begin{proof}
See Section \ref{sec:proofs}.
\end{proof}

{\sc Note 2: }{\it Limits in Corollary \eqref{Tuite} are equivalent to 
\begin{align*}
\lim_{n\to\infty}\left|\frac{\alpha_k(n+1)-\alpha_k(n)}{\alpha_k(n)-\alpha_k(n-1)}\right|,
\qquad 
\lim_{n\to\infty}\left|\frac{\alpha_k(n)}{\alpha_k(n-1)}\right|,
\end{align*}
if the latter ones exist. The numbers $R_k$, $\tilde{R}_k$ represent the radius of convergence of Maclaurin series of $(1-s)G\alpha_k(s)$ and $G\alpha_k(s)$ respectively.
}

Let us turn to the limit $\lim_{n\to\infty}a_n$. Eq. \eqref{alpha_over_b_with_initial} for $s\in\mathbb{C}$ such that $B(s)-s^m\neq0$, $m\in\mathbb{N}$, implies
\begin{align}\label{alpha_k}
(1-s)G\alpha_k(s)=\frac{1-s}{B(s)-s^m}\sum_{n=k}^{m-1}b_{n-k}s^n,\text{ for every } k=0,\,1,\,\ldots,\,m-1.
\end{align}
Due to 
\begin{align*}
\lim_{s\to1^-}(1-s)G\alpha_k(s)&=
\lim_{s\to1^-}\left(\alpha_k(0)+\sum_{n=1}^{\infty}(\alpha_k(n)-\alpha_k(n-1))s^n\right)\\
&=\alpha_k(0)+\alpha_k(1)-\alpha_k(0)+\alpha_k(2)-\alpha_k(1)+\ldots,
\end{align*}
which is valid for every $k=0,\,1,\,\ldots,\,m-1$, eq. \eqref{alpha_k} suggests the conditions when the finite limit of every sequence $\alpha_k(n)$ exists as $n\to\infty$. On the other hand, because of \eqref{a_n_expr}, the finite limit of $a_n$ as $n\to\infty$ and $m\geqslant2$ can exist even if the limits $\alpha_k(n)$ are infinite or do not exist. The following statement holds.

\begin{theorem}\label{thm:sum_b_to_1}
If every Maclaurin series $M_k(s)$ of $(1-s)G\alpha_k(s)$, $k=0,\,1,\,\ldots,\,m-1$ converges as $s\to1^{-}$, then 
\begin{align}\label{alpha_limit}
\lim_{s\to1^-}M_k(s)=\lim_{n\to\infty}\alpha_k(n),\text{ for all } k=0,\,1,\,\ldots,\,m-1.
\end{align}

In addition to \eqref{alpha_limit}, if $b_0+b_1+\ldots=1$ and $m\neq\sum_{j=1}^{\infty}j\,b_j$, then
\begin{align}\label{alpha_limits_express}
\lim_{n\to\infty}\alpha_k(n)=\frac{\sum\limits_{j=k}^{m-1}b_{j-k}}{m-\sum\limits_{j=1}^{\infty}j b_j},\,k=0,\,1,\,\ldots,\,m-1.
\end{align}
\end{theorem}

\begin{proof}
See Section \ref{sec:proofs}.
\end{proof}

If every sequence $\alpha_0(n),\,\alpha_0(n),\,\ldots,\,\alpha_{m-1}(n)$ has the finite limit as $n\to\infty$, then, according to \eqref{a_n_expr},
\begin{align}\label{a_alpha_limit}
\lim_{n\to\infty}a_n=a_0\lim_{n\to\infty}\alpha_0(n)
+a_1\lim_{n\to\infty}\alpha_1(n)+\ldots+a_{m-1}\lim_{n\to\infty}\alpha_{m-1}(n).
\end{align}

Eq. \eqref{a_alpha_limit} connects $\lim_{n\to\infty}a_n$ to the linear combination of $a_0,\,a_1,\,\ldots,\,a_{m-1}$. In other words, if eq. \eqref{a_alpha_limit} is valid, we may choose $a_0,\,a_1,\,\ldots,\,a_{m-1}$ freely and, having a fixed sequence ${\pmb b}$, know the limit of $a_n$ in the recurrence \eqref{eq:seq_v1} as $n\to\infty$. On the other hand, if eq. \eqref{a_alpha_limit} is valid, we may choose $\lim_{n\to\infty}a_n$ and ask what initial values $a_0,\,a_1,\,\ldots,\,a_{m-1}$ have to be chosen that ${\pmb a}$ from\eqref{eq:seq_v1} tends towards where we want. Thus, the last thing we conduct in this analysis is the expressions of every single term out of $a_0,\,a_1,\,\ldots,\,a_{m-1}$ via $\lim_{n\to\infty}a_n$. In general, there are no unique values $a_0,\,a_1,\,\ldots,\,a_{m-1}$ that determine $\lim_{n\to\infty}a_n$.

Let $\alpha_1,\,\alpha_2,\,\ldots$ be the roots of $B(s)-s^m$, $m\in\mathbb{N}$. The fundamental theorem of algebra guarantees that there are at least $m$ such roots counted with their multiplicities. Our last theorem and corollary express every single term out of $a_0,\,a_1,\,\ldots,\,a_{m-1}$ via $\lim_{n\to\infty}a_n$. It is obvious that if $m=1$, then $\lim_{n\to\infty}a_n=a_0\lim_{n\to\infty}\alpha_0(n)$.

\begin{theorem}\label{thm:expresions_general}
Let $m\geqslant2$. Assume that the sequence of numbers ${\pmb b}$ is such that $b_0\neq0$, $b_0+b_1+\ldots=1$ and $m\neq\sum_{j=1}^{\infty}jb_j$. Let $\alpha_1,\,\alpha_2,\,\ldots,\,\alpha_{m-1}\neq1$ be the simple roots of $B(s)-s^m$ and say that $a_0,\,a_1,\,\ldots,\,a_{m-1}$ are such that the finite limit $\lim_{n\to\infty}a_n$ exists and denote
\begin{align*}
\lim_{s\to\alpha_j}\frac{B(s)-s^m}{1-s}\sum_{k=0}^{m-1}a_kM_k(s)&=:L_j,\, j=1,\,2,\,\ldots,\,m-1,\\
\left(m-\sum_{j=1}^{\infty}jb_j\right)\lim_{n\to\infty}a_n&=:L.
\end{align*}
Here $M_k(s)$ are Maclaurin series of $(1-s)G\alpha_k(s)$, $k=0,\,1,\,\ldots,\,\ldots,\,m-1$.
Then
\begin{align}\label{system_general}
\hspace{-0.6cm}\begin{pmatrix}
\sum\limits_{n=0}^{m-1}b_n\alpha_1^n&\sum\limits_{n=1}^{m-1}b_{n-1}\alpha_1^n&\ldots&b_0\alpha_1^{m-2}+b_1\alpha_1^{m-1}&b_0\alpha_1^{m-1}\\
\sum\limits_{n=0}^{m-1}b_n\alpha_2^n&\sum\limits_{n=1}^{m-1}b_{n-1}\alpha_2^n&\ldots&b_0\alpha_2^{m-2}+b_1\alpha_2^{m-1}&b_0\alpha_2^{m-1}\\
\vdots&\vdots&\ddots&\vdots&\vdots\\
\sum\limits_{n=0}^{m-1}b_n\alpha_{m-1}^n&\sum\limits_{n=1}^{m-1}b_{n-1}\alpha_{m-1}^n&\ldots&b_0\alpha_{m-1}^{m-2}+b_1\alpha_{m-1}^{m-1}&b_0\alpha_{m-1}^{m-1}\\
\sum\limits_{n=0}^{m-1}b_n&\sum\limits_{n=1}^{m-1}b_{n-1}&\ldots&b_0+b_1&b_0
\end{pmatrix}
\begin{pmatrix}
a_0\\
a_1\\
\vdots\\
a_{m-2}\\
a_{m-1}
\end{pmatrix}
=\begin{pmatrix}
L_1\\
L_2\\
\vdots\\
L_{m-1}\\
L
\end{pmatrix}
\end{align}
and the system \eqref{system_general} has the unique solution $(a_0,\,a_1,\,\ldots,\,a_{m-1})$.
\end{theorem}

\begin{proof}
See Section \ref{sec:proofs}.
\end{proof}

Theorem \ref{thm:expresions_general} implies the following corollary.
\newpage
\begin{cor}\label{cor:expressions}
If the conditions of Theorem \ref{thm:expresions_general} are satisfied and, in addition, $L_1=L_2=\ldots=L_{m-1}=0$, then
\begin{align}
&a_0=\frac{L}{b_0}\prod_{j=1}^{m-1}\frac{\alpha_j}{\alpha_j-1},\label{sol_0}\\
&a_1=-\frac{b_0+b_1}{b_0}a_0\label{sol_1}+\frac{L}{b_0}\prod_{j=1}^{m-1}\frac{1}{\alpha_j-1}\left(\prod_{j=1}^{m-1}\alpha_j-\sum\limits_{1\leqslant j_1<\ldots<j_{m-2}\leqslant m-1}
\alpha_{j_1}\cdots\alpha_{j_{m-2}}\right),\\
&a_2=-\frac{b_0+b_1}{b_0}a_1-\frac{b_0+b_1+b_2}{b_0}a_0+\frac{L}{b_0}\prod_{j=1}^{m-1}\frac{1}{\alpha_j-1}\label{sol_2}\\
&\times\left(\prod_{j=1}^{m-1}\alpha_j-\sum_{1\leqslant j_1<\ldots<j_{m-2}\leqslant m-1}
\alpha_{j_1}\cdots\alpha_{j_{m-2}}
+\sum_{1\leqslant j_1<\ldots<j_{m-3}\leqslant m-1}
\alpha_{j_1}\cdots\alpha_{j_{m-3}}\right),\nonumber\\
&\,\,\vdots\nonumber\\
&a_{m-1}=-\frac{1}{b_0}\sum_{i=0}^{m-2}\left(\sum_{k=0}^{m-1-i}b_k\right)a_i
+\frac{L}{b_0}\prod_{j=1}^{m-1}\frac{1}{\alpha_j-1}
\times\Bigg(\prod_{j=1}^{m-1}\alpha_j-\label{sol_D-1}\\
&\sum_{1\leqslant j_1<\ldots<j_{m-2}\leqslant m-1}
\alpha_{j_1}\cdots\alpha_{j_{m-2}}
+\sum_{1\leqslant j_1<\ldots<j_{m-3}\leqslant m-1}
\alpha_{j_1}\cdots\alpha_{j_{m-3}}+\ldots+(-1)^{m+1}\Bigg),\nonumber
\end{align}
\end{cor}

\begin{proof}
See Section \ref{sec:proofs}.
\end{proof}

The above provided relations among sequences ${\pmb a}$, ${\pmb b}$, and $\alpha_0(n)$, $ \alpha_1(n)$, ... $\alpha_{m-1}(n)$, $n\in\mathbb{N}_0$, $m\in\mathbb{N}$ imply various interesting and picturesque examples that can be further studied. We provide many such examples in Section \ref{sec:examples}. However, we would like to highlight a few of them in this section as well. In Propositions \ref{ex_zeta_1}-\ref{ex_zeta_2}, we set $m=1$ and provide such sequences ${\pmb b}$ whose terms add up to one, $\alpha_0(n),\,n\in\mathbb{N}_0$ gives the partial sums of the Riemann zeta function and the ''expectation'' of nonnegative integers $b_1+2b_2+3b_3+\ldots$ computes the values of the Riemann zeta function $\zeta(s)$.  

\begin{proposition}\label{ex_zeta_1}
Let $\zeta(s)=\sum_{n=1}^{\infty}n^{-s}$, $\Re s>1$ be the Riemann zeta function. Let $a\in\mathbb{C}$ be such that $\Re a>1$ and define the sequence
\begin{align}\label{b_for_zeta_1}
b_0=1,\,b_1=-\frac{1}{2^a},\,b_n=\frac{1}{n^{a}}-\sum_{j=0}^{n-1}\frac{b_j}{(n+1-j)^{a}},\,n=2,\,3,\,\ldots
\end{align}
Then $b_0+b_1+\ldots=1$, $\alpha_0(n)=\sum_{j=1}^{n+1}j^{-a},\,n\in\mathbb{N}_0$ and
\begin{align}\label{zeta_av_1}
1-\frac{1}{\zeta(a)}=\sum_{j=1}^{\infty}j b_j=-\frac{1}{2^a}+\left(\frac{1}{2^a}-\frac{1}{3^a}+\frac{1}{4^a}\right)2+\left(\frac{1}{3^a}-\frac{2}{4^a}+\frac{2}{6^a}-\frac{1}{8^a}\right)3+\ldots,
\end{align}
where $a\in\mathbb{C}$ is such that $\Re a>1$.

Moreover, if $\Re a>1$ and we define the sequence
\begin{align}\label{b_for_zeta_2}
\tilde{b}_0=1,\,\tilde{b}_1=\frac{1}{2^a},\,\tilde{b}_n=\frac{\mu(n)}{n^{a}}-\sum_{j=0}^{n-1}\frac{\mu(n+1-j)\tilde{b}_j}{(n+1-j)^a},\,n=2,\,3,\,\ldots,
\end{align}
where $\mu(n)$ is the Möbius function (see \cite[p. 91]{Edwards}), then $\tilde{b}_0+\tilde{b}_1+\ldots=1$, $\tilde{\alpha}_0(n)=\sum_{j=1}^{n+1}\mu(j)j^{-a}$, $n\in\mathbb{N}_0$ and
\begin{align}\label{zeta_av_2}
\zeta(a)&=1-\sum_{j=1}^{\infty}j\tilde{b}_j\nonumber \\
&=1-\frac{1}{2^a}-\left(\frac{1}{3^a}-\frac{1}{2^a}+\frac{1}{4^a}\right)2-
\left(\frac{2}{6^a}-\frac{1}{4^a}-\frac{1}{8^a}-\frac{1}{3^a}\right)3-\ldots
,\,\Re \alpha>1.
\end{align}
\end{proposition}

\begin{proof}
See Section \ref{sec:proofs}.
\end{proof}

Let us denote the formal power series
\begin{align}\label{zeta_auxiliar}
f(s,\,a):=(1-s)/\sum_{n=0}^{\infty}\frac{1}{2^n}\sum_{j=0}^{n}\binom{n}{j}\frac{(-1)^j}{(j+1)^a}s^n.
\end{align}
The next proposition computes the values of the Riemann zeta function in the entire complex plane except for some of its points.

\begin{proposition}\label{ex_zeta_2}
Let $\zeta(s)$ be the Riemann zeta function defined for $\Re(s)>1$ and ana\-ly\-ti\-cal\-ly continued to the whole complex plane except the point $s=1$. Let $a\in\mathbb{C}$ and define the sequence
\begin{align*}
&\hat{b}_0=1,\,\hat{b}_1=2^{-a-1}-2^{-1},\\
&\hat{b}_n=\frac{1}{2^{n-1}}\sum_{j=0}^{n-1}\binom{n-1}{j}\frac{(-1)^j}{(j+1)^{a}}-\sum_{l=0}^{n-1}\frac{\hat{b}_l}{2^{n-l}}\sum_{j=0}^{n-l}\binom{n-l}{j}\frac{(-1)^j}{(j+1)^{a}},\,n=2,\,3,\,\ldots
\end{align*}

If $a\in\mathbb{C}$ is such that $f(s,\,a)$ from \eqref{zeta_auxiliar} is analytic in $|s|<1$ and $\lim_{s\to1^-}$f(s,\,a)=0, then $\hat{b}_0+\hat{b}_1+\ldots=1$,
\begin{align*}
\hat{\alpha}_0(n)&=\sum_{k=0}^{n}\frac{1}{2^k}\sum_{j=0}^{k}\binom{k}{j}\frac{(-1)^j}{(j+1)^a},\,n\in\mathbb{N}_0,
\end{align*}
and
\begin{align*}
\zeta(a)=\frac{1}{2-2^{2-a}}\cdot\frac{1}{1-\sum_{j=1}^{\infty}j \hat{b}_j},\, a\in\mathbb{C}\setminus\left\{1+\frac{2\pi i}{\log 2}n\right\},\,n\in\mathbb{Z}.
\end{align*}
\end{proposition}

\begin{proof}
See Section \ref{sec:proofs}.
\end{proof}

Proposition \ref{ex_zeta_2} implies 
\begin{align*}
&1-\frac{1}{(2-2^{2-a})\zeta(a)}=\sum_{j=1}^{\infty}j \hat{b}_j\\
&=\frac{1}{2^{a+1}}-\frac{1}{2}+\left(\frac{1}{2}-\frac{1}{2^{a+1}}+\frac{1}{4^{a+1}}-\frac{1}{4\cdot3^a}\right)2 +\left(\frac{1}{8\cdot3^a}-\frac{1}{4\cdot6^a}+\frac{1}{8^{a+1}}\right)3+\ldots,
\end{align*}
where $a\in\mathbb{C}$ is such that $(2-2^{2-a})\zeta(a)\neq0$. 

Moreover, Proposition \ref{ex_zeta_2} implies the following equivalent condition for the Riemann hypothesis.

\begin{cor}\label{cor_RH}
The Riemann hypothesis is true iff $\hat{b}_1+2\hat{b}_2+\ldots$ converges in $1/2<\Re a<1$.
\end{cor}

{\sc Note 3:} The convergence of $\hat{b}_1+2\hat{b}_2+\ldots$ can be associated with the vanishing order of modulus of coefficients $\hat{b}_0,\, \hat{b}_1,\,\ldots$ For instance, if $|\hat{b}_n|\sim n^{-2-\varepsilon}$, where $\varepsilon>0$ and $n\to\infty$, then this series converges.

In our last Proposition \ref{proposition_pi} below, we set $m=1$ again and provide such sequences ${\pmb b}$ that $\alpha_0(n),\,n\in\mathbb{N}_0$ corresponds to the partial sums that approximate $\pi$ and $e$ respectively.

\begin{proposition}\label{proposition_pi}
Let 
\begin{align*}
\bar{b}_0=1,\,\bar{b}_1=\frac{1}{3},\,\bar{b}_n=\frac{(-1)^{n+1}}{2n-1}+\sum_{j=1}^{n}\frac{(-1)^{j+1}\,\bar{b}_{n-j}}{2j+1},\,n=2,\,3,\,\ldots
\end{align*}
Then $\bar{b}_0+\bar{b}_1+\ldots=1$,
\begin{align*}
\bar{\alpha}_0(n)=\sum_{j=1}^{n+1}\frac{(-1)^{j+1}}{2j-1},\,n\in\mathbb{N}_0,
\end{align*}
and
\begin{align*}
1-\frac{4}{\pi}=\sum_{j=1}^{\infty}j\bar{b}_j=\frac{1}{3}-\frac{19}{45}\cdot2+\frac{128}{945}\cdot3-\frac{1088}{14175}\cdot4+\frac{4864}{93555}\cdot5-\ldots
\end{align*}
Moreover, if we define
\begin{align*}
\dot{b}_0=1,\,\dot{b}_1=-1,\,\dot{b}_n=\frac{1}{(n-1)!}-\sum_{j=1}^{n}\frac{\dot{b}_{n-j}}{j!},\,n=2,\,3,\,\ldots,
\end{align*}
then $\dot{b}_0+\dot{b}_1+\dots=1$, $\dot{\alpha}_0(n)=\sum_{j=0}^{n}1/j!$, $n\in\mathbb{N}_0$ and
\begin{align*}
1-\frac{1}{e}=\sum_{j=1}^{\infty}j\dot{b}_j=-1+\frac{3}{2}\cdot2-\frac{2}{3}\cdot3+\frac{5}{24}\cdot4-\frac{1}{20}\cdot5+\frac{7}{720}\cdot6-\frac{1}{630}\cdot7+\ldots
\end{align*}

\end{proposition}

\begin{proof}
See Section \ref{sec:proofs}.
\end{proof}

As mentioned, sometimes the chosen sequences ${\pmb b}$ determine ${\pmb a}$ or $\alpha_0(n)$, $ \alpha_1(n)$, ... $\alpha_{m-1}(n)$, $n\in\mathbb{N}_0$, $m\in\mathbb{N}$ such that these sequences represent some strictly structured attractive patterns in plane or space. Here are some examples depicted with \cite{Mathematica} when $m=1$:  

\begin{figure}[H]
    \begin{minipage}[t]{.3\textwidth}
        \centering
        \includegraphics[width=\textwidth]{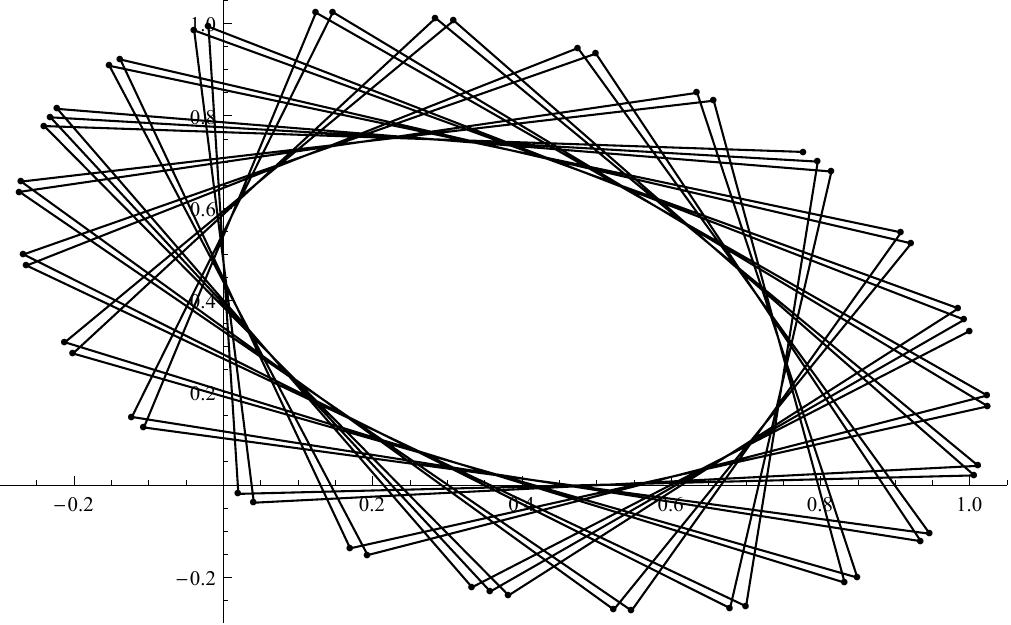}
        \caption*{(a) $N=50$}
    \end{minipage}
    \hfill
    \begin{minipage}[t]{.3\textwidth}
        \centering
        \includegraphics[width=\textwidth]{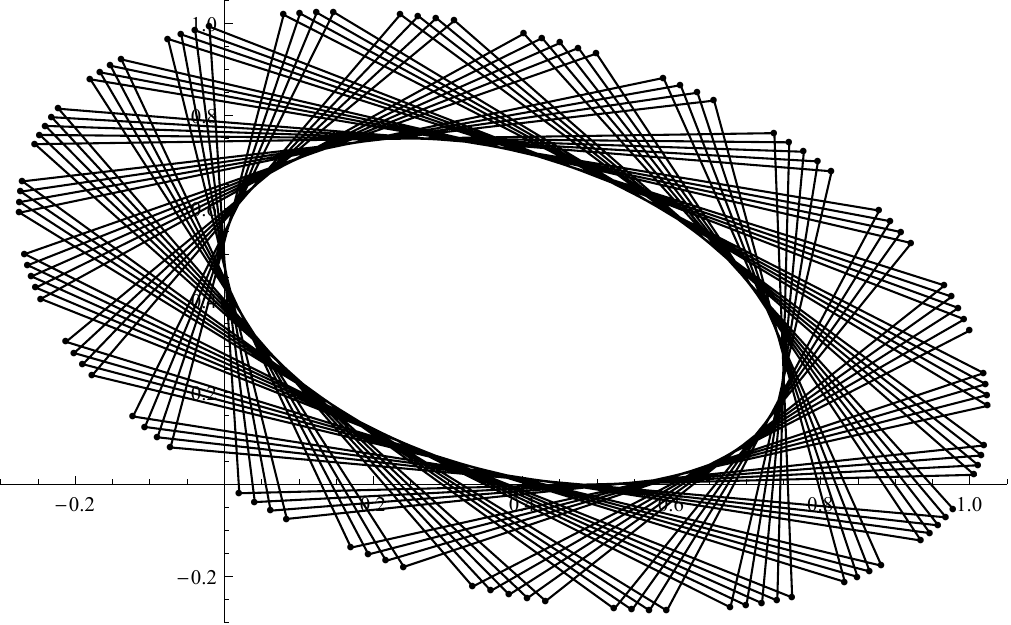}
        \caption*{(b) $N=100$}
    \end{minipage}
    \hfill
    \begin{minipage}[t]{.3\textwidth}
        \centering
        \includegraphics[width=\textwidth]{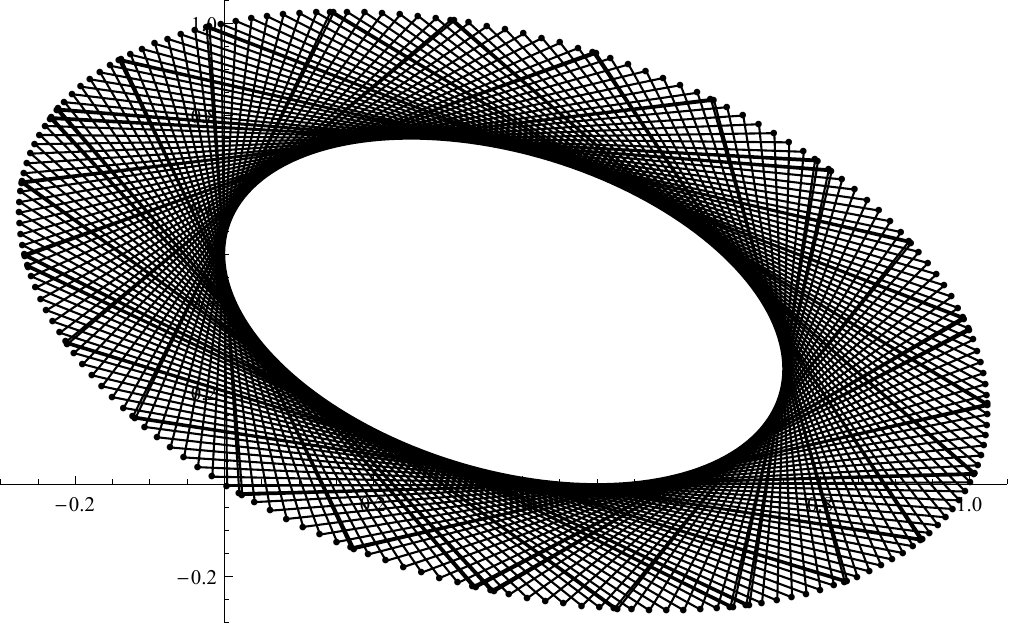}
        \caption*{(c) $N=200$}
    \end{minipage}  
      \caption{${\pmb b}=\{-3,\,2,\,-1,\,3,\,0,\,0,\,\ldots,\,0,\,\ldots\}$ and the pairs of $(\alpha_0(n),\,\alpha_0(n+1))$ in plane connected with lines when $n$ varies from $0$ to $N$.}\label{F1}
\end{figure}

\begin{figure}[H]
    \begin{minipage}[t]{.3\textwidth}
        \centering
        \includegraphics[width=\textwidth]{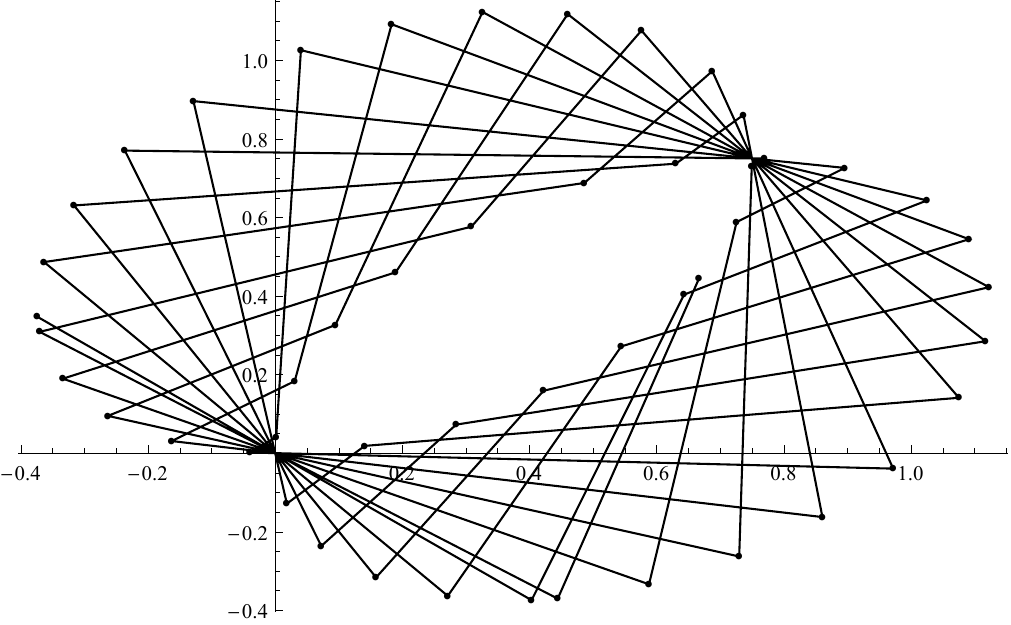}
        \caption*{(a) $N=50$}
    \end{minipage}
    \hfill
    \begin{minipage}[t]{.3\textwidth}
        \centering
        \includegraphics[width=\textwidth]{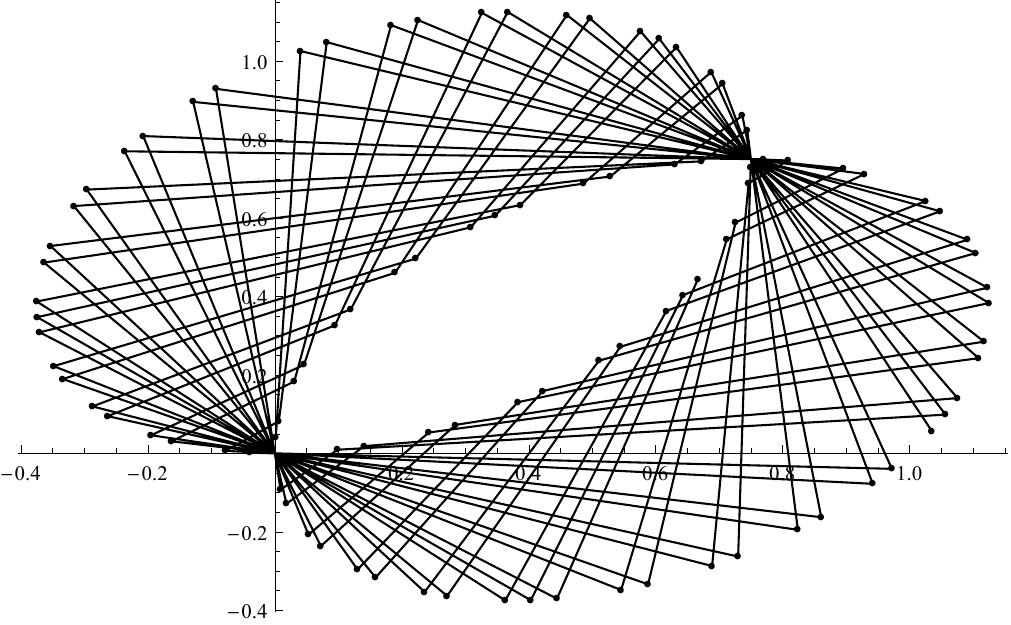}
        \caption*{(b) $N=100$}
    \end{minipage}
    \hfill
    \begin{minipage}[t]{.3\textwidth}
        \centering
        \includegraphics[width=\textwidth]{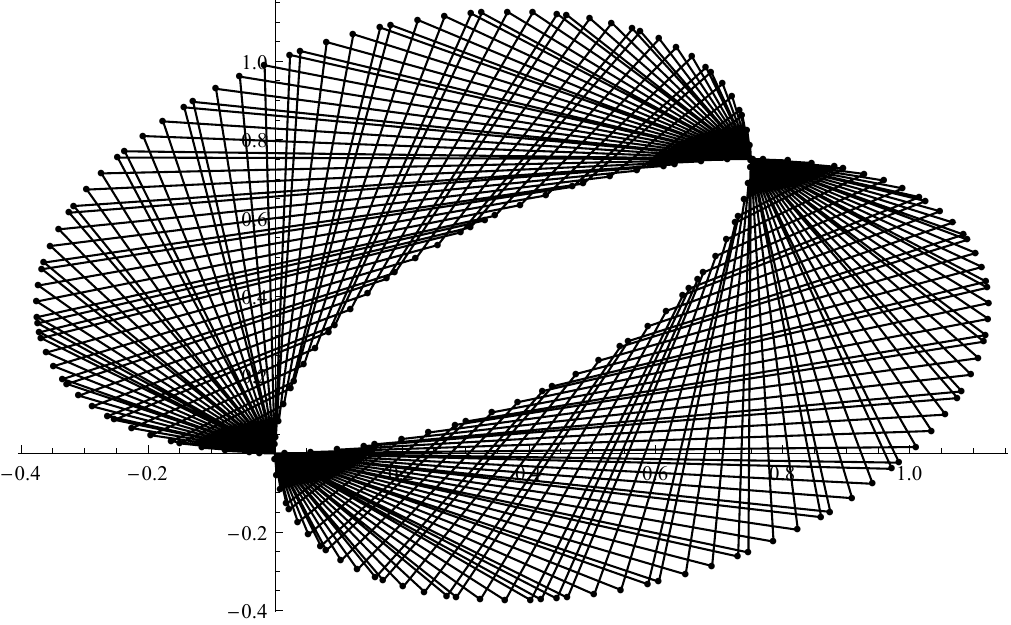}
        \caption*{(c) $N=200$}
    \end{minipage}  
       \caption{${\pmb b}=\{3,\,-1,\,0,\,2,\,-3,\,0,\,0,\,\ldots,\,0,\,\ldots\}$ and the pairs of $(\alpha_0(n),\,\alpha_0(n+1))$ in plane connected with lines when $n$ varies from $0$ to $N$.}
\end{figure}

\begin{figure}[H]
    \begin{minipage}[t]{.3\textwidth}
        \centering
        \includegraphics[width=\textwidth]{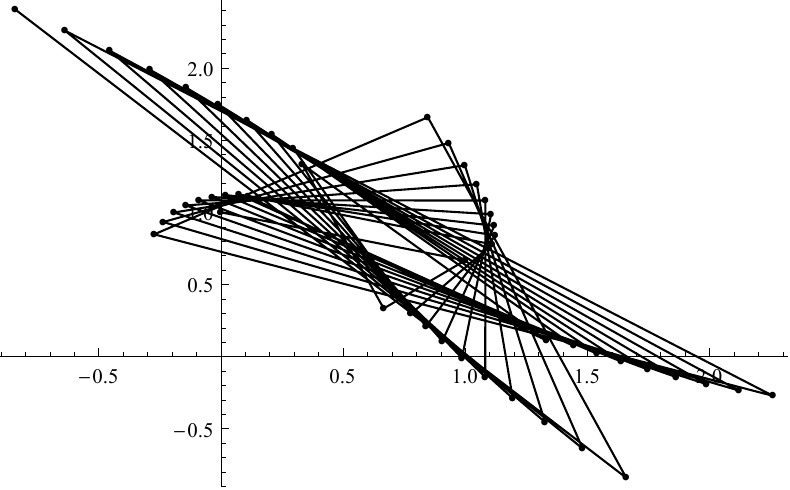}
        \caption*{(a) $N=50$}
    \end{minipage}
    \hfill
    \begin{minipage}[t]{.3\textwidth}
        \centering
        \includegraphics[width=\textwidth]{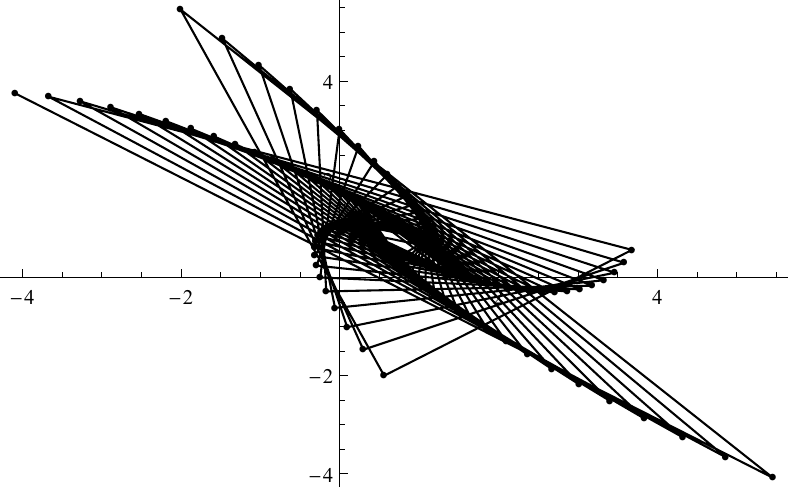}
        \caption*{(b) $N=100$}
    \end{minipage}
    \hfill
    \begin{minipage}[t]{.3\textwidth}
        \centering
        \includegraphics[width=\textwidth]{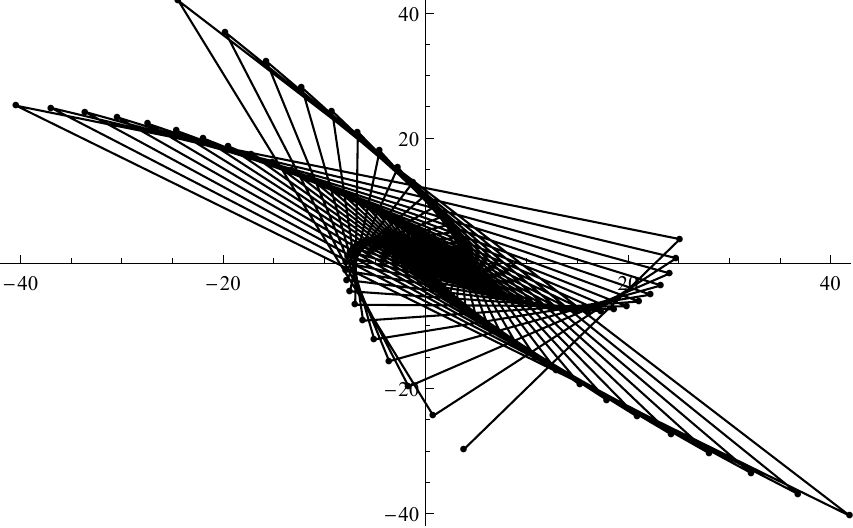}
        \caption*{(c) $N=200$}
    \end{minipage}  
       \caption{${\pmb b}=\{3,\,1,\,-3,\,-2,\,2,\,0,\,0,\,\ldots,\,0,\,\ldots\}$ and the pairs of $(\alpha_0(n),\,\alpha_0(n+1))$ in plane connected with lines when $n$ varies from $0$ to $N$.}
\end{figure}

\begin{figure}[H]
    \begin{minipage}[t]{.3\textwidth}
        \centering
        \includegraphics[width=\textwidth]{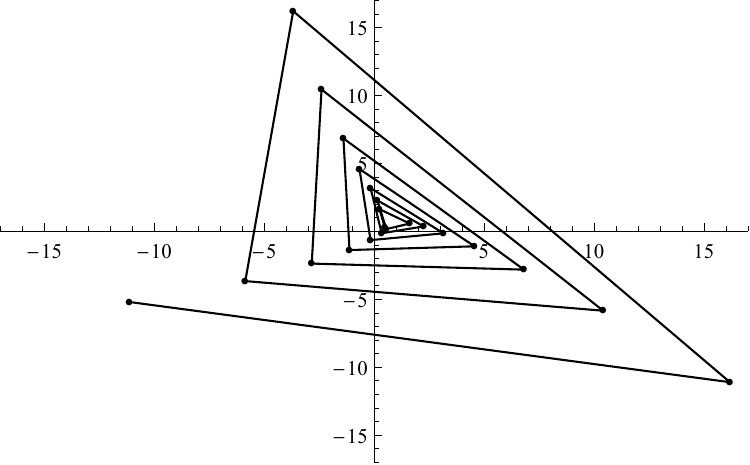}
        \caption*{(a) $N=25$}
    \end{minipage}
    \hfill
    \begin{minipage}[t]{.3\textwidth}
        \centering
        \includegraphics[width=\textwidth]{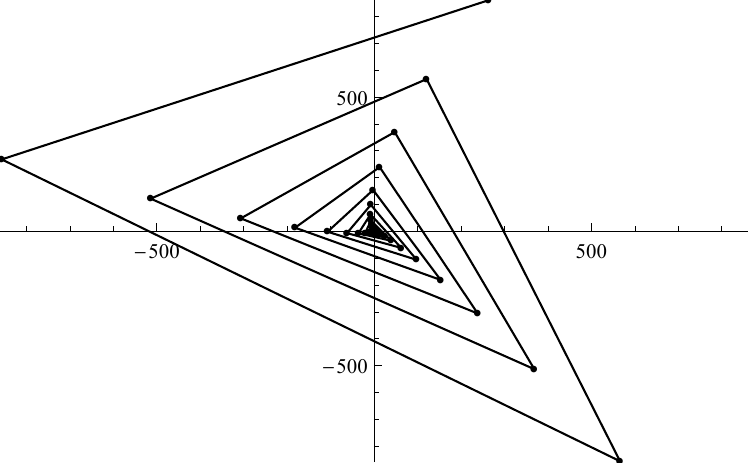}
        \caption*{(b) $N=50$}
    \end{minipage}
    \hfill
    \begin{minipage}[t]{.3\textwidth}
        \centering
        \includegraphics[width=\textwidth]{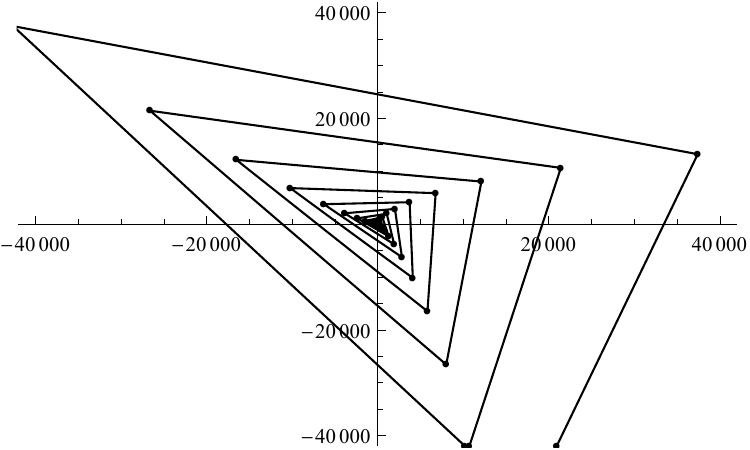}
        \caption*{(c) $N=75$}
    \end{minipage}  
       \caption{${\pmb b}=\{2,\,0,\,0,\,-3,\,2,\,0,\,0,\,\ldots,\,0,\,\ldots\}$ and the pairs of $(\alpha_0(n),\,\alpha_0(n+1))$ in plane connected with lines when $n$ varies from $0$ to $N$.}
\end{figure}

\begin{figure}[H]
    \begin{minipage}[t]{.3\textwidth}
        \centering
        \includegraphics[width=\textwidth]{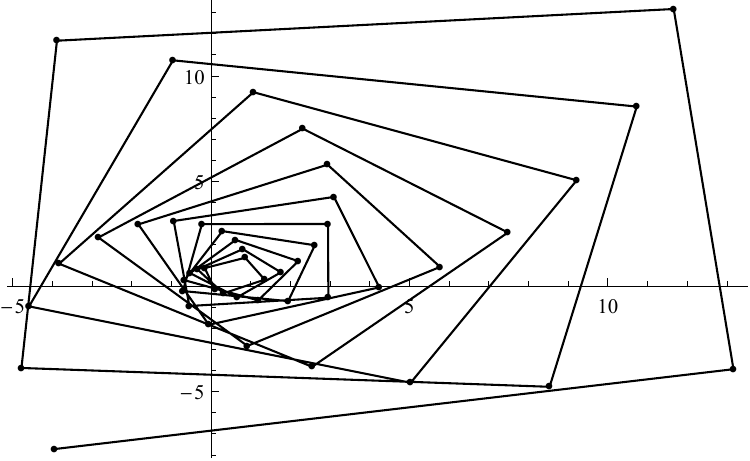}
        \caption*{(a) $N=50$}
    \end{minipage}
    \hfill
    \begin{minipage}[t]{.3\textwidth}
        \centering
        \includegraphics[width=\textwidth]{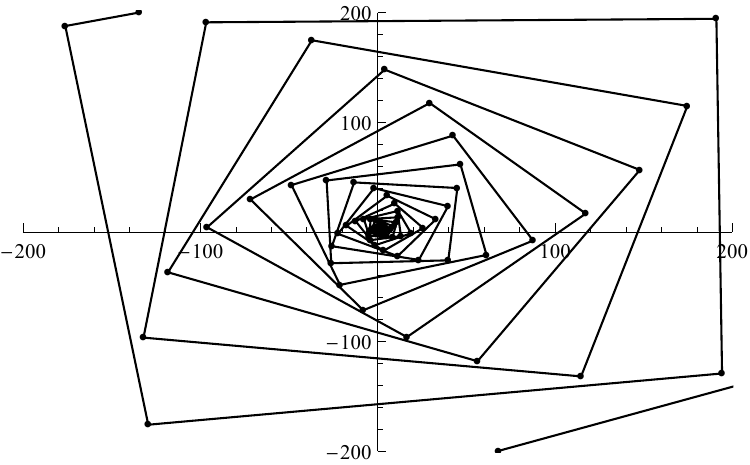}
        \caption*{(b) $N=100$}
    \end{minipage}
    \hfill
    \begin{minipage}[t]{.3\textwidth}
        \centering
        \includegraphics[width=\textwidth]{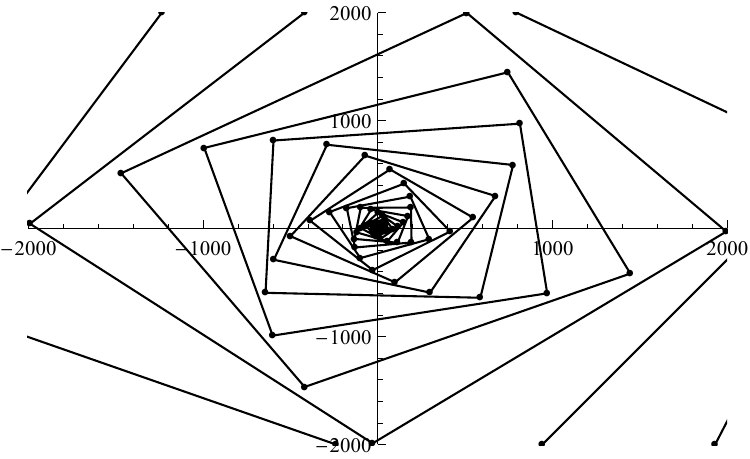}
        \caption*{(c) $N=150$}
    \end{minipage}  
       \caption{${\pmb b}=\{\sin1,\,\sin2,\,\sin3,\,\ldots\}$ and the pairs of $(\alpha_0(n),\,\alpha_0(n+1))$ (in radians) in plane connected with lines when $n$ varies from $0$ to $N$.}
\end{figure}

\begin{figure}[H]
    \begin{minipage}[t]{.3\textwidth}
        \centering
        \includegraphics[width=\textwidth]{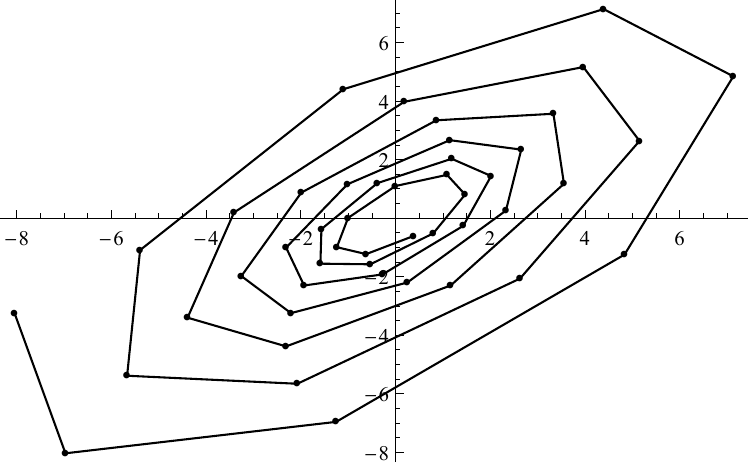}
        \caption*{(a) $N=50$}
    \end{minipage}
    \hfill
    \begin{minipage}[t]{.3\textwidth}
        \centering
        \includegraphics[width=\textwidth]{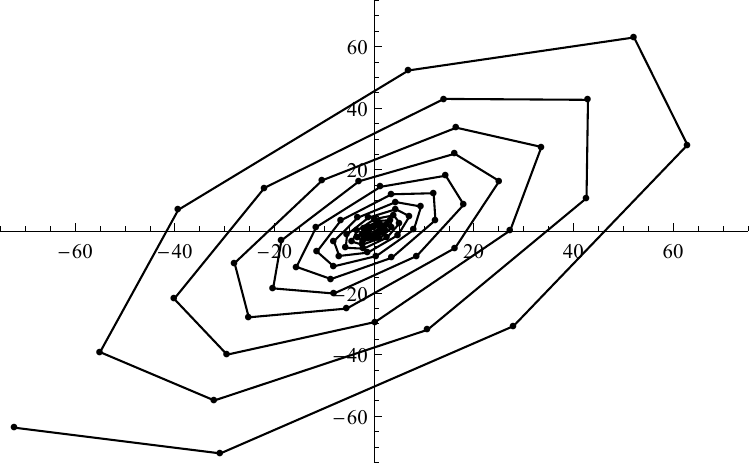}
        \caption*{(b) $N=100$}
    \end{minipage}
    \hfill
    \begin{minipage}[t]{.3\textwidth}
        \centering
        \includegraphics[width=\textwidth]{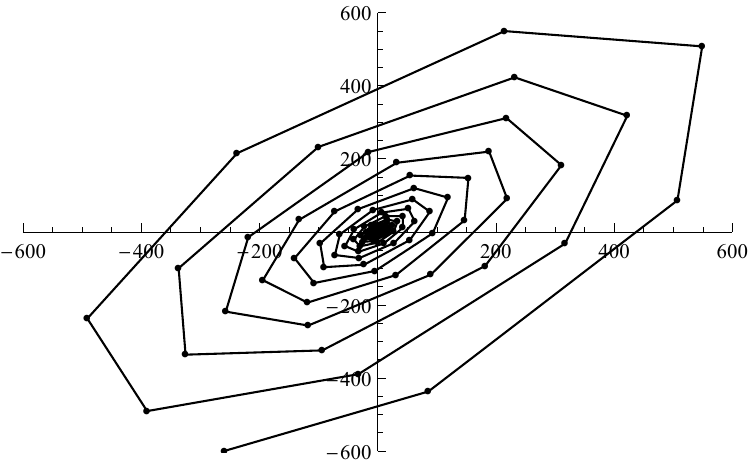}
        \caption*{(c) $N=150$}
    \end{minipage}  
       \caption{${\pmb b}=\{F(n+1)/\varphi^n,\,n\in\mathbb{N}_0\}$, where $F(n)$ denotes the Fibonnaci number and $\varphi$ is the golden ratio. We depict the pairs of $(\alpha_0(n),\,\alpha_0(n+1))$ in plane connected with lines when $n$ varies from $0$ to $N$.}
\end{figure}

\begin{figure}[H]
    \begin{minipage}[t]{.3\textwidth}
        \centering
        \includegraphics[width=\textwidth]{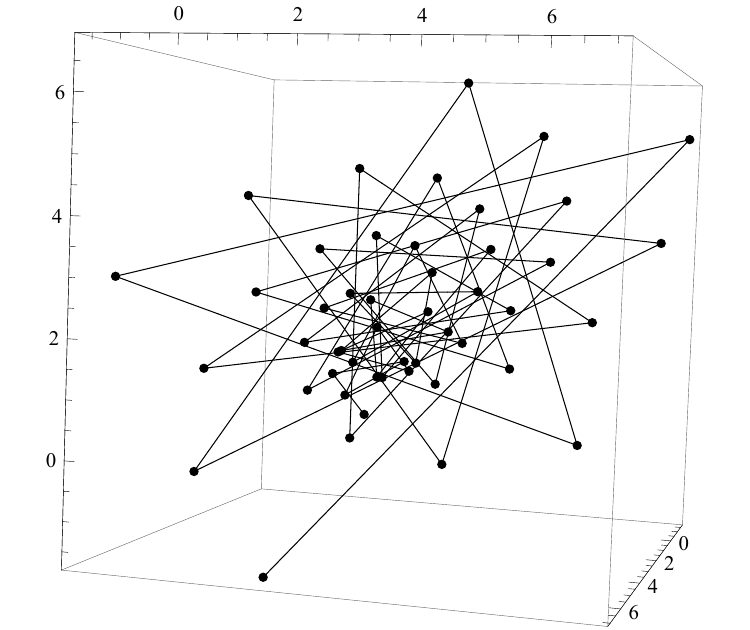}
        \caption*{(a) $N=50$}
    \end{minipage}
    \hfill
    \begin{minipage}[t]{.3\textwidth}
        \centering
        \includegraphics[width=\textwidth]{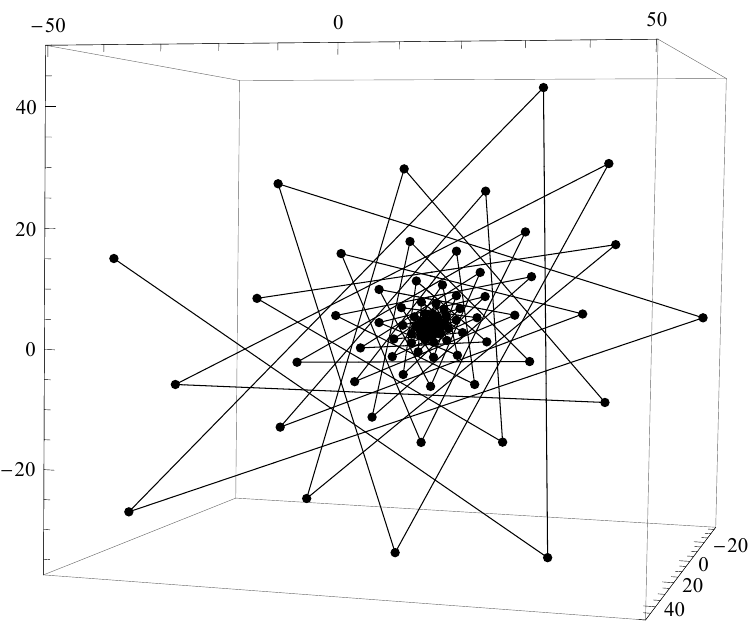}
        \caption*{(b) $N=100$}
    \end{minipage}
    \hfill
    \begin{minipage}[t]{.3\textwidth}
        \centering
        \includegraphics[width=\textwidth]{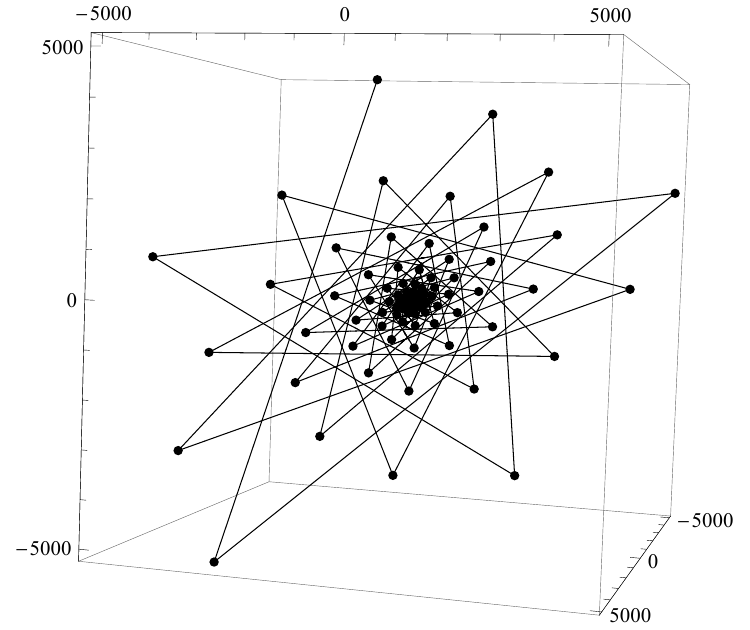}
        \caption*{(c) $N=200$}
    \end{minipage}  
\caption{${\pmb b}=\{3,\,0,\,-3,\,-2,\,3,\,0,\,0,\,\ldots,\,0,\,\ldots\}$, and the triples of $(\alpha_0(n),\,\alpha_0(n+1),\,\alpha_0(n+2))$ in space connected with lines when $n$ varies from $0$ to $N$.}\label{F7}
\end{figure}

\section{Proofs of the main results}\label{sec:proofs}

\begin{proof}[Proof of Proposition \ref{lem}]
The proof is based on mathematical induction. If $n=0,\,1,\,\ldots,\,m-1$, the equality \eqref{a_n_expr} is evident in view of Table \ref{T} realizing that it states $a_0=a_0$, $a_1=a_1$, $\ldots$, $a_{m-1}=a_{m-1}$. If $n=m$, the equality \eqref{a_n_expr} is the same as \eqref{eq:seq_v1}. If $n=m+1,\,m+2,\,\ldots$, then by \eqref{eq:seq_v1} and induction hypothesis
\begin{align*}
a_n&=\frac{1}{b_0}\left(a_{n-m}-\sum_{j=0}^{n-1}b_{n-j}a_j\right)
=\frac{1}{b_0}\left(\sum_{k=0}^{m-1}\alpha_k(n-m)a_k-\sum_{j=0}^{n-1}b_{n-j}\sum_{k=0}^{m-1}\alpha_k(j)a_k\right)\\
&=\frac{1}{b_0}\left(\alpha_0(n-m)-\sum\limits_{j=0}^{n-1}b_{n-j}\alpha_{0}(j)\right)a_0
+\frac{1}{b_0}\left(\alpha_1(n-m)-\sum\limits_{j=0}^{n-1}b_{n-j}\alpha_{1}(j)\right)a_1\\
&+\ldots+\frac{1}{b_0}\left(\alpha_{m-1}(n-m)-\sum\limits_{j=0}^{n-1}b_{n-j}\alpha_{m-1}(j)\right)a_{m-1}
=\sum_{k=0}^{m-1}\alpha_{k}(n)a_k.
\end{align*}
\end{proof}

\begin{proof}[Proof of Theorem \ref{thm}]
Eq. \eqref{eq:seq} implies 
\begin{align*}
a_{n-m}=\sum_{j=0}^{n}b_{n-j}a_j\quad &\Rightarrow \quad
\sum_{n=m}^{\infty}a_{n-m}\,s^{n-m}=\sum_{n=m}^{\infty}\left(\sum_{j=0}^{n}b_{n-j}a_j\right)s^{n-m}\\
&\Rightarrow \quad
\sum_{j=0}^{m-1}a_j\sum_{n=j}^{m-1}b_{n-j}s^n=A(s)(B(s)-s^m),\,m\in\mathbb{N}.
\end{align*}

This proves \eqref{a_over_b}. Eq. \eqref{a_n_expr} of Proposition \ref{lem} yields \eqref{a_over_alpha} immediately.

Arguing the same as deriving \eqref{a_over_b}, for every recurrence $\alpha_0(j)$, $\alpha_1(j)$, $\ldots$, $\alpha_{m-1}(j)$, $j\in\mathbb{N}_0$, $m\in\mathbb{N}$, defined in \eqref{rec_for_expres}, we obtain
\begin{align}\label{alpha_over_b}
\sum_{j=0}^{m-1}\alpha_k(j)\sum_{n=j}^{m-1}b_{n-j}s^n=G\alpha_k(s)(B(s)-s^m),\,m\in\mathbb{N},
\end{align}
where $k=0,\,1,\,\ldots,\,m-1$. By inserting the initial values of recurrences $\alpha_0(j)$, $\alpha_1(j)$, $\ldots$, $\alpha_{m-1}(j)$ from Table \ref{T} into \eqref{alpha_over_b}, we obtain 
\begin{align*}
\sum_{n=k}^{m-1}b_{n-k}s^n=G\alpha_k(s)(B(s)-s^m),\, m\in\mathbb{N},
\end{align*}
where $k=0,\,1,\,\ldots,\,m-1$. This proves \eqref{alpha_over_b_with_initial}. Eq. \eqref{a_diff} is implied by Proposition \ref{lem} immediately computing $a_n-a_{n-1}$, $n\in\mathbb{N}$ from there.

For $s\in\mathbb{C}$ such that $B(s)\neq s^m$, \eqref{alpha_over_b_with_initial} yields
\begin{align}\label{alpha_over_b_with_initial_s-1}
(1-s)\cdot\frac{\sum\limits_{n=k}^{m-1}b_{n-k}s^n}{B(s)-s^m}=(1-s)G\alpha_k(s)=
\alpha_k(0)+\sum_{n=1}^{\infty}(\alpha_k(n)-\alpha_k(n-1))s^n,
\end{align}
for every $k=0,\,1,\,\ldots,\,m-1$. Equating the coefficients at the powers of \eqref{alpha_over_b_with_initial_s-1}, we get \eqref{alpha_diff}. Eq. \eqref{alpha_ne_diff} is implied by \eqref{alpha_over_b_with_initial} immediately.  
\end{proof}

\begin{proof}[Proof of Corollary \ref{Tuite}]
The functions $(1-s)G\alpha_k(s)$ and $G\alpha_k(s)$ are defined at $s=0$ and analytic in some circles of positive radius all centered at the origin. Therefore, $(1-s)G\alpha_k(s)$ and $G\alpha_k(s)$ coincide with their Maclaurin series in the described circles. Thus, $R_k$ or $\tilde{R}_k$ are also the radius of convergence of the Maclaurin series of $(1-s)G\alpha_k(s)$ or $G\alpha_k(s)$, see \cite{complex_analysis}.
\end{proof}

\begin{proof}[Proof of Theorem \ref{thm:sum_b_to_1}]
If $s\to1^-$, the right-hand side of \eqref{alpha_over_b_with_initial_s-1} is $\lim_{n\to\infty}\alpha_k(n)$, for every $k=0,\,1,\,\ldots,\,m-1$. It remains to set conditions for the left-hand side of \eqref{alpha_over_b_with_initial_s-1} such that the finite limits $\lim_{n\to\infty}\alpha_k(n)$ exist. This proves \eqref{alpha_limit}. 

If $b_0+b_1+\ldots=1$, then
\begin{align*}
&B(s)-s^m=b_0+\sum_{n=1}^{\infty}b_n s^n-s^m=1-s^m-\sum_{n=1}^{\infty}b_n(1-s^n)\\
&=(1-s)\sum_{j=0}^{m-1}s^j-(1-s)\sum_{n=1}^{\infty}b_n\sum_{j=0}^{n-1}s^j
=(1-s)\left(\sum_{j=0}^{m-1}s^j-\sum_{j=0}^{\infty}s^j\sum_{n=j+1}^{\infty}b_n\right).
\end{align*}

By inserting the above expression into \eqref{alpha_over_b_with_initial}, for every $k=0,\,1,\,\ldots,\,m-1$ we obtain
\begin{align}\label{series_before_limit_m}
\frac{\sum\limits_{n=k}^{m-1}b_{n-k}\,s^{n}}{\sum\limits_{j=0}^{m-1}s^j-\sum\limits_{j=0}^{\infty}s^j\sum\limits_{n=j+1}^{\infty}b_n}=(1-s)G\alpha_k(s).
\end{align}

Thus, \eqref{alpha_limits_express} follows by letting $s\to1^-$ in the both sides of \eqref{series_before_limit_m} under the requirement that the finite limits
\begin{align*}
\lim_{s\to1^-}(1-s)G\alpha_k(s)=\lim_{n\to\infty}\alpha_k(n).
\end{align*}
exist.
\end{proof}

\begin{proof}[Proof of Theorem \ref{thm:expresions_general}]
The system \eqref{system_general} is implied by \eqref{a_over_b}. Indeed, if $m\in\{2,\,3,\,\ldots\}$ and $\alpha_1,\,\alpha_2,\,\ldots,\,\alpha_{m-1}\neq1$ are the simple roots of $B(s)-s^m$, then the first $m-1$ out of $m$ equations of \eqref{system_general} are implied by letting $s\to\alpha_k$, $k=1,\,2,\,\ldots,\,m-1$ in the both sides of \eqref{a_over_b}. Notice that the Maclaurin series of $(1-s)A(s)$ is $\sum_{k=0}^{m-1}a_kM_k(s)$, where $M_k(s)$ are the Maclaurin series of $(1-s)G\alpha_k(s)$. The last equation of \eqref{system_general} follows from \eqref{a_over_b} as well. Indeed, if we let $s\to1^-$ in both sides of \eqref{a_over_b}, then
\begin{align*}
\sum_{k=0}^{m-1}a_k\sum_{n=k}^{m-1}b_{n-k}=\lim_{s\to1^-}(1-s)A(s)\frac{B(s)-s^m}{1-s}=\left(m-\sum_{j=1}^{\infty}jb_j\right)\lim_{n\to\infty}a_n.
\end{align*}
The system \eqref{system_general} has a unique solution because its matrix determinant is
\begin{align}\label{determinant}
\frac{b_0^m}{(-1)^{m+1}}\prod\limits_{j=1}^{m-1}(\alpha_j-1)\prod\limits_{1\leqslant i <j \leqslant m-1}(\alpha_j-\alpha_i),\,m\in\mathbb{N}\setminus\{1\},
\end{align}
see \cite[Lem. 4.2]{A1} for more details on how this determinant is computed. The conditions when the provided determinant \eqref{determinant} $\neq0$ are evident.
\end{proof}

\begin{proof}[Proof of Corollary \ref{cor:expressions}]

If $L_1=L_2=\ldots=L_{m-1}=0$ and $L\neq0$ in system \eqref{system_general}, then, to get its solution, it suffices to know the minors of the last line of the system's matrix. These minors are expressible as a combination of the roots $\alpha_1,\,\alpha_2,\,\ldots,\,\alpha_{m-1}$ similarly as determinant \eqref{determinant}. For more details, see \cite[Proof of  Thm. 3.3]{A1}.
\end{proof}

\begin{proof}[Proof of Proposition \ref{ex_zeta_1}]
Let $s\in\mathbb{C}$ and denote
\begin{align*}
\Phi(s,\,a,\,2):=\sum_{n=0}^{\infty}\frac{s^n}{(n+2)^a},\,\Re a>1,\,|s|<1.
\end{align*}

Then, the recurrence relation \eqref{b_for_zeta_1}, implies
\begin{align*}
&\sum_{n=2}^{\infty}s^nb_n=\sum_{n=2}^{\infty}\frac{s^n}{n^a}-\sum_{n=2}^{\infty}\left(\sum_{j=0}^{n-1}\frac{b_j}{(n+1-j)^a}\right)s^n\quad \Rightarrow\\
&B(s)-1+\frac{s}{2^a}=s^2\Phi(s,\,a,\,2)-sB(s)\Phi(s,\,a,\,2)+\frac{s}{2^a}\quad \Rightarrow\\
&B(s)=\frac{1+s^2\Phi(s,\,a,\,2)}{1+s\Phi(s,\,a,\,2)}=\frac{1-s}{1+s\Phi(s,\,a,\,2)}+s,\,|s|<1.
\end{align*}
The last expression implies $\lim_{s\to1^-}B(s)=b_0+b_1+\ldots=1$. If we choose, $m=1$, then, according to \eqref{alpha_over_b_with_initial},
\begin{align*}
G\alpha_0(s)=\frac{1}{B(s)-s}=\frac{1+s\Phi(s,\,a,\,2)}{1-s},\,|s|<1,
\end{align*}
whose Maclaurin series yields
\begin{align*}
\alpha_0(n)=\sum_{j=0}^{n}\frac{1}{(j+1)^a},\,n\in\mathbb{N}_0.
\end{align*}
Indeed, these coefficients $\alpha_0(n),\,n\in\mathbb{N}_0$ ''sit'' at the powers of $s$ of the series
\begin{align*}
\left(1+\sum_{n=0}^{\infty}\frac{s^{n+1}}{(n+2)^a}\right)(1+s+s^2+\ldots),\,|s|<1,\,\Re a>1.
\end{align*}

According to Theorem \ref{thm:sum_b_to_1}
\begin{align*}
\lim_{n\to\infty}\alpha_0(n)=\zeta(a)=\frac{1}{1-\sum_{j=1}^{\infty}jb_j},\,\Re a>1
\end{align*}
and the formula \eqref{zeta_av_1} is proved.

Formula \eqref{zeta_av_2} follows analogously. Indeed, sequence \eqref{b_for_zeta_2} for some $s\in\mathbb{C}$ implies 
\begin{align*}
\tilde{B}(s)=\sum_{n=0}^{\infty}\tilde{b}_ns^n=\frac{1-s}{1+s\sum\limits_{n=0}^{\infty}\frac{\mu(n+2)}{(n+2)^a}s^n}+s.
\end{align*}
Then, choosing $m=1$, from \eqref{alpha_over_b_with_initial} we get
\begin{align*}
\tilde{G}\alpha_0(s)=\frac{1}{\tilde{B}(s)-s}=\frac{1+s\sum\limits_{n=0}^{\infty}\frac{\mu(n+2)}{(n+2)^a}s^n}{1-s}
\quad \Rightarrow \quad \tilde{\alpha}_0(n)=\sum_{j=0}^{n}\frac{\mu(j+1)}{(j+1)^a},\,n\in\mathbb{N}_0,
\end{align*}
and, by Theorem \ref{thm:sum_b_to_1},
\begin{align*}
\lim_{n\to\infty}\tilde{\alpha}_0(n)=\frac{1}{\zeta(a)}=\frac{1}{1-\sum_{j=1}^{\infty}\tilde{b}_j j},\,\Re a >1.
\end{align*}
\end{proof}

\begin{proof}[Proof of Proposition \ref{ex_zeta_2}]
The proof is nearly identical to that of Proposition \ref{ex_zeta_1}. After long and careful derivation (we omit details), we obtain the formal power series
\begin{align*}
\hat{B}(s)&=\sum_{n=0}^{\infty}s^n \hat{b}_n=\frac{1-s}{\sum\limits_{n=0}^{\infty}\frac{1}{2^n}\sum\limits_{j=0}^{n}\binom{n}{j}\frac{(-1)^j}{(j+1)^a}s^n}+s.
\end{align*}
Under conditions of Proposition \ref{ex_zeta_2}, $\lim_{s\to1^-}\hat{B}(s)=\hat{b}_0+\hat{b}_1+\ldots=1$.

Choosing $m=1$, from \eqref{alpha_over_b_with_initial} we get
\begin{align*}
\hat{G}\alpha_0(s)&=\frac{1}{\hat{B}(s)-s}=\frac{1}{1-s}\sum\limits_{n=0}^{\infty}\frac{1}{2^n}\sum\limits_{j=0}^{n}\binom{n}{j}\frac{(-1)^j}{(j+1)^a}s^n,\,|s|<1,
\end{align*}
whose coefficients of its Maclaurin series are 
\begin{align*}
\hat{\alpha}_0(n)&=\sum_{k=0}^{n}\frac{1}{2^k}\sum_{j=0}^{k}\binom{k}{j}\frac{(-1)^j}{(j+1)^a},\,n\in\mathbb{N}_0.
\end{align*}

Due to the Hasse formula for the Riemann zeta function \cite{Hasse} 
\begin{align}\label{Hasse}
\zeta(a)=\frac{1}{1-2^{1-a}}\sum_{k=0}^{\infty}\frac{1}{2^k}\sum_{j=0}^{k}\binom{k}{j}\frac{(-1)^j}{(j+1)^a},\,a\neq1+\frac{2\pi i}{\log2}n,\,n\in\mathbb{Z},
\end{align}
and Theorem \ref{thm:sum_b_to_1}
\begin{align*}
\lim_{n\to\infty}\hat{\alpha}_0(n)&=(2-2^{2-a})\zeta(a)=\frac{1}{1-\sum_{j=1}^{\infty}j\hat{b}_j}.
\end{align*}
\end{proof}

\begin{proof}[Proof of Proposition \ref{proposition_pi}]
Arguing the same as proving Propositions \ref{ex_zeta_1} and \ref{ex_zeta_2}, we obtain
\begin{align*}
\bar{B}(s)=\sum_{n=1}^{\infty}\bar{b}_ns^n=\frac{\sqrt{s}(1-s)}{\arctan\sqrt{s}}+s,\,|s|<1.
\end{align*}
Notice that $s/\arctan s\to1$ as $s\to0$. Choosing $m=1$, from \eqref{alpha_over_b_with_initial} we get
\begin{align*}
\bar{G}\alpha_0(s)=\frac{\arctan{\sqrt{s}}}{\sqrt{s}(1-s)},\, |s|<1,
\end{align*}
whose coefficients of its Maclaurin series are
\begin{align*}
\bar{\alpha}_0(n)=\sum_{k=1}^{n+1}\frac{(-1)^{k+1}}{2k-1},\,n\in\mathbb{N}_0.
\end{align*}

Due to the Leibniz formula for $\pi$ \cite{pi} and Theorem \ref{thm:sum_b_to_1}
\begin{align*}
\lim_{n\to\infty}\bar{\alpha}_0(n)&=\frac{\pi}{4}=\frac{1}{1-\sum_{j=1}^{\infty}j\bar{b}_j}.
\end{align*}
For the last sequence $\dot{b}_0,\,\dot{b}_1,\,\ldots$, when $|s|<1$, we have
\begin{align*}
\dot{B}(s)=\frac{1-s}{e^s}+s,\,\dot{G}\alpha_0(s)=\frac{e^s}{1-s},\,
\dot{\alpha}_0(n)=\sum_{j=0}^{n}\frac{1}{j!},\,n\in\mathbb{N}_0,
\end{align*}
and the rest is the same as previously.
\end{proof}

\section{Examples}\label{sec:examples}
Let us recall the notations ${\pmb a}:=(a_0,\,a_1,\,\ldots)$ and ${\pmb b}:=(b_0,\,b_1,\,\ldots)$. We also denote ${\pmb \alpha}_k:=\{\alpha_k(0),\,\alpha_k(1),\,\ldots\}$, $k=0,\,1,\,\ldots,\,m-1$, $m\in\mathbb{N}$. In this section, we provide more examples than those in Propositions \ref{ex_zeta_1}--\ref{proposition_pi} when chosen sequences ${\pmb b}$ determine ${\pmb a}$ or ${\pmb \alpha}_k$. All the necessary computer computations are implemented with \cite{Mathematica}. We emphasize that details, when obtaining $B(s)$ and the Maclaurin series of $B(s)$-based functions in some particular cases, are extensive and are omitted in this section. However, the omitted details are available from the authors on demand.

\subsection{Cases when the limit exists}

\begin{ex}\label{1}
Let ${\pmb b}=\{5,\,-4,\,-3,\,3,\,0,\,0,\,\ldots,\,0,\,\ldots\}$ and $m=1$. We determine the sequence ${\pmb a}$ that is related to ${\pmb b}$ as \eqref{eq:seq_v1}.
\end{ex}

Let $M_0(s)$ and $G_0(s)$ correspondingly be the Maclaurin series of $(1-s)G\alpha_0(s)$ and $G\alpha_0(s)$ whose radiuses of convergence are denoted by $R_M$ and $R_G$. Then, according to \eqref{alpha_over_b_with_initial},
\begin{align*}
(1-s)G\alpha_0(s)&=\frac{1}{1-3/5s^2},\,G\alpha_0(s)=\frac{1}{(1-s)(1-3/5s^2)} \\
M_0(s)&=\sum_{n=0}^{\infty}\left(\frac{3}{5}\right)^n s^{2n}, R_M=\sqrt{\frac{5}{3}}\approx1.29099,\\
G_0(s)&=\frac{5}{2}\sum_{n=0}^{\infty}\left(1-\left(\frac{3}{5}\right)^{ \floor*{n/2}+1}\right)s^n,\,R_G=1.
\end{align*}

Here $\floor{\cdot}$ denotes the floor function. Then, according to \eqref{alpha_ne_diff},
\begin{align*}
\alpha_0(n)=\frac{5}{2}\left(1-\left(\frac{3}{5}\right)^{ \floor*{n/2}+1}\right),\,n\in\mathbb{N}_0\quad \Rightarrow \quad 
\lim_{n\to\infty}\alpha_0(n)=\frac{5}{2},
\end{align*}
where the limit of $\alpha_0(n)$ can also be obtained by \eqref{alpha_limit} as $M_0(1)=5/2$.

In view of \eqref{a_n_expr}, $a_n=\alpha_0(n)a_0$, $n\in\mathbb{N}_0$ and $\lim_{n\to\infty}a_n=5/2\cdot a_0$. So if we choose $a_0=2/5$, then
\begin{align*}
a_n=\left(1-\left(\frac{3}{5}\right)^{ \floor*{n/2}+1}\right),\,n\in\mathbb{N}_0.
\end{align*}

\begin{ex}
Let ${\pmb b}$ be as in Example \ref{1} and say that $m=2$. We construct the sequence ${\pmb a}$ whose values are computed as in \eqref{eq:seq_v1}.
\end{ex}

Let $M_0(s)$ and $M_1(s)$ correspondingly be the Maclaurin series of $(1-s)G\alpha_0(s)$ and $(1-s)G\alpha_1(s)$. According to \eqref{alpha_over_b_with_initial}, we have that 
\begin{align*}
(1-s)G\alpha_0(s)&=\frac{5-4s}{5+s-3s^2},\,M_0(s)=\sum_{n=0}^{\infty}c_ns^n,\\
(1-s)G\alpha_1(s)&=\frac{5s}{5+s-3s^2},\,M_1(s)=\sum_{n=0}^{\infty}d_ns^n,
\end{align*}
where 
\begin{align*}
&c_0=1,\,c_1=-1,\,c_{n}=(3c_{n-2}-c_{n-1})/5,\,n=2,\,3\,\ldots,\\
&d_0=0,\,d_1=1,\,d_{n}=(3d_{n-2}-d_{n-1})/5,\,n=2,\,3\,\ldots,
\end{align*}
and the radius of convergence of both series $M_0(s)$ and $M_1(s)$ is $R_M=(-1+\sqrt{61})/6\approx1.13504$; $R_M$ is the smallest root in modulus of $1/(1-s)/G\alpha_0(s)$ and $1/(1-s)/G\alpha_1(s)$. Then, because of \eqref{alpha_diff},
\begin{align}\label{aswellas}
\begin{cases}
\alpha_0(0)=1,\,\alpha_0(n)=\alpha_0(n-1)+c_n\\
\alpha_1(0)=0,\,\alpha_1(n)=\alpha_1(n-1)+d_n
\end{cases},\,
n\in\mathbb{N}.
\end{align}

According to \eqref{alpha_limit}, $M_0(1)=1/3=\lim_{n\to\infty}\alpha_0(n)$ and $M_1(1)=5/3=\lim_{n\to\infty}\alpha_1(n)$. Due to \eqref{a_n_expr}, $a_n=a_0\alpha_0(n)+a_1\alpha_1(n)$, $n\in\mathbb{N}_0$ and
\begin{align}\label{limit_in_example}
\lim_{n\to\infty}a_n=\frac{1}{3}a_0+\frac{5}{3}a_1.
\end{align}

On the other hand, if $G_0(s)$ and $G_1(s)$ correspondingly denote the Maclaurin series of $G\alpha_0(s)$ and $G\alpha_1(s)$, then, according to \eqref{alpha_over_b_with_initial},
\begin{align*}
G\alpha_0(s)&=\frac{5-4s}{(1-s)(5+s-3s^2)},\,G_0(s)=\sum_{n=0}^{\infty}\alpha_0(n)s^n,\\
G\alpha_1(s)&=\frac{5s}{(1-s)(5+s-3s^2)},\,G_1(s)=\sum_{n=0}^{\infty}\alpha_1(n)s^n,
\end{align*}
where 
\begin{align*}
\begin{cases}
\alpha_0(0)=1,\,\alpha_0(1)=0,\,\alpha_0(2)=\frac{4}{5},\,\alpha_0(n)=\frac{4\alpha_0(n-1)+4\alpha_0(n-2)-3\alpha_0(n-3)}{5}\\
\alpha_1(0)=0,\,\alpha_1(1)=1,\,\alpha_1(2)=\frac{4}{5},\,\alpha_1(n)=\frac{4\alpha_1(n-1)+4\alpha_1(n-2)-3\alpha_1(n-3)}{5}
\end{cases},\,
n=3,\,4\,\ldots
\end{align*}
and the radius of convergence of both series $G_0(s)$ and $G_1(s)$ is $R_G=1$; $R_G$ is the smallest root in modulus of $1/G\alpha_0(s)$ and $1/G\alpha_1(s)$.
This, as well as \eqref{aswellas}, yields
\begin{align*}
{\pmb \alpha}_0&=\left\{1,\,0,\,\frac{4}{5},\,\frac{1}{25},\,\frac{84}{125},\,\frac{56}{625},\,\frac{1829}{3125},\,\frac{2136}{15625},\,\ldots\right\},\\
{\pmb \alpha}_1&=\left\{0,\,1,\,\frac{4}{5},\,\frac{36}{25},\,\frac{149}{125},\,\frac{1016}{625},\,\frac{4344}{3125},\,\frac{26521}{15625},\,\ldots\right\},
\end{align*}
and, according to \eqref{R},
\begin{align*}
\limsup_{n\to\infty}
\left|\alpha_0(n+1)-\alpha_0(n)\right|^{1/n}
=
\limsup_{n\to\infty}\left|\alpha_{1}(n+1)-\alpha_1(n)\right|^{1/n}
=\frac{1}{R_M}&=\frac{6}{-1+\sqrt{61}}\\
&\approx0.88102.
\end{align*}
Let's see what $a_0$ and $a_1$ in \eqref{limit_in_example} can be chosen to express both of them via $\lim_{n\to\infty}a_n$. The roots $\neq1$ of $B(s)-s^2$ are
\begin{align*}
B(s)=5-4s-3s^2+3s^3=s^2\quad \Rightarrow \quad \alpha_{1}:=\frac{1-\sqrt{61}}{6},\,
\alpha_2:=\frac{1+\sqrt{61}}{6}.
\end{align*}
The both series $M_0(s)$, $M_1(s)$ do not converge as $s\to\alpha_1$ or $s\to\alpha_2$.

Assuming that $a_0$ and $a_1$ are such that
\begin{align*}
\lim_{s\to\alpha_1}(a_0M_0(s)+a_1M_1(s))\frac{B(s)-s^2}{1-s}=0,
\end{align*}
by \eqref{sol_0} and \eqref{sol_1} we get
\begin{align*}
a_0&=\frac{2-B'(1)}{b_0}\cdot\frac{\alpha_1}{\alpha_1-1}\cdot \lim_{n\to\infty}a_n=\frac{11-\sqrt{61}}{10}\cdot \lim_{n\to\infty}a_n,\\
a_1&=-\frac{b_0+b_1}{b_0}\cdot a_0+\frac{2-B'(1)}{b_0}\cdot\lim_{n\to\infty}a_n=\frac{19+\sqrt{61}}{50}\cdot\lim_{n\to\infty}a_n,
\end{align*}
where $B'(\ldots)$ denotes derivative.

If we require $a_0$ and $a_1$ be such that
\begin{align*}
\lim_{s\to\alpha_2}(a_0M_0(s)+a_1M_1(s))\frac{B(s)-s^2}{1-s}=0,
\end{align*}
\bigskip
then, by \eqref{sol_0} and \eqref{sol_1}, we get
\begin{align*}
a_0&=\frac{2-B'(1)}{b_0}\cdot\frac{\alpha_1}{\alpha_1-1}\cdot \lim_{n\to\infty}a_n=\frac{11+\sqrt{61}}{10}\cdot \lim_{n\to\infty}a_n,\\
a_1&=-\frac{b_0+b_1}{b_0}\cdot a_0+\frac{2-B'(1)}{b_0}\cdot\lim_{n\to\infty}a_n=\frac{19-\sqrt{61}}{50}\cdot\lim_{n\to\infty}a_n.
\end{align*}

\begin{ex}
Let
\begin{align*}
b_n = (-1)^{n+1}\frac{\pi^{2n}}{(2n)!}\left(1+\frac{i \pi}{2n+1}\right),\,n\in\mathbb{N}_0.
\end{align*}
where the sum $b_0+b_1+\ldots=1$ represents the well-known Euler's identity $e^{\pi i}=-1$. We chose $m=1$ and compute the sequence ${\pmb \alpha}_0$.
\end{ex}
In this example,
\begin{align*}
& B(s)= -\cos(\pi\sqrt{s})-\frac{i\sin(\pi\sqrt{s})}{\sqrt{s}}, \\
&(1-s)G\alpha_0(s)= \frac{b_0(1-s)}{B(s)-s} = \frac{(1+i\pi)(1-s)}{\cos(\pi\sqrt{s})+\frac{i\sin(\pi\sqrt{s})}{\sqrt{s}}+s}, \\
&M_0(s) = 1 +\frac{12i-6\pi-3i\pi^2+\pi^3}{6(\pi-i)}s\\
&+\frac{-720-360i\pi+540\pi^2+360i\pi^3-135\pi^4-42i\pi^5+7\pi^6}{360(\pi-i)^2}s^2 + \dots
=\sum\limits_{n=0}^{\infty}c_ns^n,
\end{align*}
where
\begin{align*}
c_0 = 1,\, c_1 = \frac{12i-6\pi-3i\pi^2+\pi^3}{6(\pi-i)},\, c_n = \frac{1}{b_0}\left(c_{n-1}-\sum_{k=1}^{n}b_k\,c_{n-k}\right),\,n=2,\,3,\,\ldots
\end{align*}
The radius of convergence of $M_0(s)$ is $R_M\approx 1.26751$; this is the modulus of $\approx 0.71491 + 1.04665i$, which is the smallest root in modulus of $1/(1-s)/G\alpha_0(s)$. Then, according to \eqref{alpha_limit},
\begin{align*}
\lim_{s\to1}M_0(s)=\lim_{n\to\infty}\alpha_0(n)=\frac{2(\pi-i)}{\pi + 2i}\approx 1.13480 - 1.35906i.
\end{align*}
According to \eqref{alpha_diff} and \eqref{R}
\begin{align*}
\alpha_0(0)&=1,\,\alpha_0(n)=\alpha_0(n-1)+c_n,\,n=1,\,2,\,\ldots,\\
{\pmb \alpha}_0&=\left\{1,\,\frac{6i-3i\pi^2+\pi^3}{6(\pi-i)},\,\frac{-360+360\pi^2+120i\pi^3-75\pi^4-42i\pi^5+7\pi^6}{360(\pi-i)^2} ,\, \ldots \right\} \\
&\approx\{1, \, 1.85560 - 0.66183i, \, 1.89477 - 1.51773i, \,\ldots\},\\
&\limsup_{n\to\infty}\left|\alpha_0(n+1)-\alpha_0(n)\right|^{1/n}=\frac{1}{R_M}\approx0.78895.
\end{align*}
In this example, let us observe one beautiful constant
\begin{align*}
\left|\lim_{n\to\infty}\alpha_0(n)\right|=\left|\frac{2(\pi-i)}{\pi+2i}\right|=2\sqrt{\frac{1+\pi^2}{4+\pi^2}}\approx1.771.
\end{align*}

\begin{ex}
Let
\begin{align*}
b_n = \frac{3\sqrt{3}}{2\pi}\frac{\Gamma^2(n+1)}{\Gamma(2n+2)}, \, n\in\mathbb{N}_0.
\end{align*}
We chose $m=1$ and compute the sequence ${\pmb\alpha}_0$.
\end{ex}
In this example,
\begin{align*}
B(s)&=\frac{6\sqrt{3}\arcsin(\sqrt{s}/2)}{\pi\sqrt{4s - s^2}},\\
(1-s)G\alpha_0(s)&=\frac{b_0(1-s)}{B(s)-s} = \frac{3\sqrt{3}(1-s)\sqrt{4s-s^2}}{12\sqrt{3}\arcsin\left(\sqrt{s}/2\right)-2\pi s\sqrt{4s-s^2}},
\end{align*}
whose Maclaurin series is
\begin{align*}
M_0(s)&= 1 + \left(\frac{2\pi}{3\sqrt{3}}-\frac{7}{6}\right)s+\left(\frac{29}{180}-\frac{8\sqrt{3}\pi}{27}+\frac{4\pi^2}{27}\right)s^2 \\ &+
\left(\frac{37}{7560}+\frac{7\pi}{30\sqrt{3}}-\frac{2\pi^2}{9} + \frac{8\pi^3}{81\sqrt{3}}\right)s^3 + \dots = \sum\limits_{n=0}^{\infty}c_ns^n,
\end{align*}
where
\begin{align*}
c_0 = 1,\,
c_n &=\frac{(-1)^n}{4^n}\left(\binom{1/2}{n}+4\binom{1/2}{n-1}\right)\\
&-\sum_{k=1}^{n}\frac{(2k)!\, c_{n-k}}{16^{k}(k!)^2(2k+1)}
+\frac{2\pi}{3\sqrt{3}}\sum_{k=0}^{n-1}\frac{(-1)^k}{4^{k}}\binom{1/2}{k}c_{n-1-k},\,n\in\mathbb{N}.
\end{align*}

Here
\begin{align*}
\binom{1/2}{0}:=1,\,\binom{1/2}{n}=\frac{(-1)^{n-1}}{n\cdot 2^{2n-1}}\binom{2n-2}{n-1},\,n\in\mathbb{N}.
\end{align*}

The radius of convergence of $M_0(s)$ is $R_M\approx3.63609$ (this is the smallest root in modulus of $1/(1-s)/G\alpha_0(s)$), the series converges as $s\to1$, and
\begin{align*}
\lim_{s\to1}M_0(s)=\lim_{n\to\infty}\alpha_0(n)=\frac{27}{-18+8\sqrt{3}\pi}\approx1.05753.
\end{align*}

Let us show that in this example $\alpha_0(n)<\alpha_0(n+1)$ for all $n\in\mathbb{N}_0$. It can be seen that $b_n>0$ for all $n\in\mathbb{N}_0$. Also, it can be shown that $\lim_{s\to1^-}B(s)=b_0+b_1+\ldots=1$. Then
\begin{align*}
1=\alpha_0(0)<1+\frac{b_2+b_3+\ldots}{b_0}=\frac{1-b_1}{b_0}=\alpha_0(1)
\end{align*}
and, under the induction hypothesis and due to \eqref{rec_for_expres},
\begin{align*}
\alpha_0(n)>\alpha_0(n)(b_0+b_1+\ldots+b_n)>
b_n\alpha_0(0)+b_{n-1}\alpha_0(1)+\ldots+b_0\alpha_0(n)=\alpha_0(n-1),
\end{align*}
for all $n\in\mathbb{N}$.

Thus, according to \eqref{alpha_diff}, \eqref{R}, {\sc Note 2} with the fact that $\alpha_0(n+1)>\alpha_0(n)$ for all $n\in\mathbb{N}_0$, we have
\begin{align*}
\alpha_0(0)&=1,\,\alpha_0(n)=\alpha_0(n-1)+c_n,\,n=1,\,2,\,\ldots,\\
{\pmb \alpha}_0&=\left\{1,\,\frac{2\pi}{3\sqrt{3}}-\frac{1}{6},\,-\frac{1}{80}-\frac{2\pi}{9\sqrt{3}}+\frac{4\pi^2}{27},\,-\frac{1}{1512}+\frac{\pi}{90\sqrt{3}}-\frac{2\pi^2}{27}+\frac{8\pi^3}{81\sqrt{3}},\,\ldots\right\}\\
&\approx\{1,\,1.04253, \, 1.05354, \, 1.05646, \,\ldots\},\\
&\lim_{n\to\infty}\frac{\alpha_0(n+1)-\alpha_0(n)}{\alpha_0(n)-\alpha_0(n-1)}=\frac{1}{R_M}\approx0.27502.
\end{align*}

\subsection{Cases when the finite limit does not exist}

\begin{ex}
Let ${\pmb b}=\{2,\,1,\,-2,\,0,\,0,\,\ldots\}$ and $m=1$. We show that the limit $\lim_{n\to\infty}a_n$ does not exist. 
\end{ex}

In  this example, according to \eqref{alpha_over_b_with_initial}, 
\begin{align*}
(1-s)G\alpha_0(s)=\frac{b_0\,(1-s)}{B(s)-s}=\frac{1}{1+s}.
\end{align*}
As the Maclaurin series of $1/(1+s)$ is $1-s+s^2-s^3+\ldots$, it is evident that this series does not converge as $s\to1$. Then, according to \eqref{alpha_diff},
\begin{align*}
\alpha_0(0)=1,\,
\begin{cases}
\alpha_0(n)-\alpha_0(n-1)=1\text{ if $n$ is even}\\
\alpha_0(n)-\alpha_0(n-1)=-1\text{ if $n$ is odd}
\end{cases},\,n\in\mathbb{N}.
\end{align*}
what yields 
\begin{align*}
{\pmb \alpha_0}=\{1,\,0,\,1,\,0,\,\ldots\}\qquad\text{or} \qquad {\pmb a}=\{a_0,\,0,\,a_0,\,0,\,\ldots\}
\end{align*}
with an arbitrarily chosen number $a_0$.

\begin{ex}
Let $m=1$ and ${\pmb b}$ represent the Fibonacci-geometric probability distribution (see \cite[\href{https://oeis.org/A302922}{A302922}]{oeis}, \cite{Fibonacci_Geometric}), i.e. 
\begin{align*}
b_n=\frac{F(n+1)}{2^{n+2}},\,n\in\mathbb{N}_0,
\end{align*}
where $F(0)=0$, $F(1)=1$, $F(n)=F(n-1)+F(n-2),\,n\geqslant 2$ is the Fibonacci sequence. We show that $\lim_{n\to\infty}a_n=\infty$ and provide pattern to compute ${\pmb \alpha}_0$ or ${\pmb a}$.
\end{ex}

In this example 
\begin{align*}
B(s)=\frac{1}{4-s(2 + s)},\,(1-s)G\alpha_0(s)=\frac{b_0\,(1-s)}{B(s)-s}=\frac{-4 + s (2 + s)}{4 (-1 + s (3 + s))},
\end{align*}
and the Maclaurin series of $(1-s)G\alpha_0(s)$ is
\begin{align*}
M_0(s)=1+\frac{1}{4}\sum_{n=1}^{\infty}c_ns^n,
\end{align*}
where
\begin{align*}
c_1=10,\,c_2=33,\,c_n=3\,c_{n-1}+c_{n-2},\,n=3,\,4,\,\ldots
\end{align*}
The series $M_0(s)$ does not converge as $s\to1$; its radius of convergence is $R_M=(\sqrt{13}-3)/2\approx0.30278$, where $R_M$ is the smallest root in modulus of $1/(1-s)/G\alpha_0(s)$. Thus, by \eqref{alpha_diff}, $\alpha_0(0)=1$, $\alpha_0(n)=\alpha_0(n-1)+c_{n}/4,\,n\in\mathbb{N}$, and, according to \eqref{R_1},
\begin{align*}
\lim_{n\to\infty}\frac{\alpha_0(n)}{\alpha_0(n-1)}=\frac{1}{R_M}=\frac{2}{\sqrt{13}-3}\approx3.30278.
\end{align*}

In conclusion, $\lim_{n\to\infty}\alpha_0(n)=\infty$, $\lim_{n\to\infty}a_n=\infty$ and
\begin{align*}
4\,{\pmb \alpha}_0=\{4,\,14,\,47,\,156,\,516,\,1705,\,5632,\,18602,\,\ldots\}
\end{align*}
or $a_n=a_{n-1}+a_0\,c_n/4,\,n\in\mathbb{N}$ with arbitrarily chosen initial value $a_0$. The provided sequence $4\,{\pmb \alpha}_0$ matches the one \cite[\href{https://oeis.org/A082574}{A082574}]{oeis} omiting its firs term.

Let us switch to the probabilistic context for a moment. The sequence ${\pmb b}$ in this example can be interpreted as a probability mass function of some non-negative and integer-valued random variable $X$, i.e., $b_n=\mathbb{P}(X=n),\,n\in\mathbb{N}_0$. Then the random variable $X+2$ describes that, for example, $2,\,3,\,\ldots$ coin tosses are required to get two heads in a row. If we define 
\begin{align*}
\tilde{a}_n=\mathbb{P}(\max\{X_1,\,X_1+X_2,\,\ldots\}\leqslant n),\,n\in\mathbb{N}_0,
\end{align*}
where $X_1,\,X_2,\,\ldots$ are independent copies of $X$, then these sequences $\tilde{{\pmb a}}$ and ${\pmb b}$ are related as in \eqref{eq:seq}. It is known (see \cite[Thm. 1 (i)]{AA}) that $\tilde{a}_n=0$ for all $n\in\mathbb{N}_0$ if $\mathbb{E}X>1$. As $\mathbb{E}X=\sum_{j=1}^{\infty}jb_j=4$ here, we conclude that 
\begin{align*}
1-\tilde{a}_n=\mathbb{P}\left(\{X_1>n\}\cup\{X_1+X_2>n\}\cup\ldots\right)=1 \text{ for all $n\in\mathbb{N}$}.
\end{align*}
In other words, the random walk $\{X_1,\,X_1+X_2,\,\ldots\}$ valued on $\mathbb{N}_0$ with probabilities determined by ${\pmb b}$ hits any threshold $n$ (no matter how large) at least once with probability one.

\begin{ex}
Let
\begin{align*}
b_n = \frac{C_n}{2^{2n+1}},\,n\in\mathbb{N}_0,   
\end{align*}
where
\begin{align*}
C_n = \frac{1}{n+1}\binom{2n}{n},\,n\in\mathbb{N}_0
\end{align*}
are Catalan numbers \cite[\href{https://oeis.org/A000108}{A000108}]{oeis}. We compute the sequence ${\pmb\alpha}_0$.
\end{ex}

In this example
\begin{align*}
B(s)=\frac{1-\sqrt{1-s}}{s},\,(1-s)G\alpha_0(s)=\frac{b_0\,(1-s)}{B(s)-s} 
=\frac{1}{2}\frac{s\sqrt{1-s}}{(1+s)\sqrt{1-s}-1},
\end{align*}
and the Maclaurin series of $(1-s)G\alpha_0(s)$ is
\begin{align*}
M_0(s)=1+\frac{3}{4}s+\frac{19}{16}s^2+\frac{61}{32}s^3+\frac{787}{256}s^4+\frac{2543}{512}s^5+\ldots=
\sum_{n=0}^{\infty}c_ns^n,
\end{align*}
where
\begin{align*}
c_0=1,\,c_n=\sum_{k=1}^{n}\frac{(2k+1)!!}{2^k(k+1)!}c_{n-k},\,n=1,\,2,\,\ldots,
\end{align*}
and the double factorial function is defined as follows
\begin{align*}
n!!=
\begin{cases}
n\cdot(n-2)\cdots5\cdot3,\,n>0 \text{ odd},\\
n\cdot(n-2)\cdots4\cdot2,\,n>0 \text{ even},\\
1,\,n=-1,\,0.
\end{cases}
\end{align*}

The radius of convergence of $M_0(s)$ is $R_M=1/\varphi=(\sqrt{5}-1)/2\approx0.61803$, where $\varphi$ denotes the golden ratio, and $R_M$ is the unique root of $1/(1-s)/G\alpha_0(s)$. Thus, the series $M_0(s)$ does not converge as $s\to1$.
According to \eqref{alpha_diff} and \eqref{R_1}:
\begin{align*}
&\alpha_0(0)=1,\,\alpha_0(n)=\alpha_0(n-1)+c_n,\,n=1,\,2,\,\ldots,\\
&{\pmb \alpha}_0=\left\{1,\,1\frac{3}{2^2},\,2\frac{15}{2^4},\,4\frac{27}{2^5},\,7\frac{235}{2^8},\,12\frac{453}{2^9},\,20\frac{1875}{2^{11}},\,33\frac{3707}{2^{12}},\,\ldots\right\},\\
&\lim_{n\to\infty}\alpha_0(n)=\infty,\\
&\lim_{n\to\infty}\frac{\alpha_0(n)}{\alpha_0(n-1)}=\varphi\approx 1.61803.
\end{align*}

\begin{ex}
Let
\begin{align*}
b_n = F(n+1)\varphi^{-n},\,n\in\mathbb{N}_0
\end{align*}
where $F(n+1)$ are the Fibonacci numbers ($F(0)=0$, $F(1)=1$, $F(n)=F(n-1)+F(n-2),\,n\geqslant2$) and $\varphi$ denotes the golden ratio. We assume $m=1$ and compute the sequence ${\pmb\alpha}_0$.
\end{ex}
In this example
\begin{align*}
&B(s)=\frac{\varphi^2}{\varphi^2-\varphi s - s^2},\,
(1-s)G\alpha_0(s)=\frac{b_0(1-s)}{B(s)-s}=\frac{(s-1)(s^2+\varphi s-\varphi^2)}{\varphi^2-s\varphi^2+\varphi s^2+s^3},\\
&\lim_{n\to\infty}b_n=\frac{1}{\sqrt{5}}\approx0.447214,\,\cite[p. 62]{penguin}.
\end{align*}
Let $M_0(s)$ and $G_0(s)$ be the corresponding Maclaurin series of $(1-s)G\alpha_0(s)$ and $G\alpha_0(s)$. Then
\begin{align*}
M_0(s)&= 1-\frac{s}{\varphi}-s^2-\frac{s^3}{\varphi}+\frac{s^4}{\varphi^3}+s^5+\frac{-9+5\sqrt{5}}{2}s^6+\ldots=\sum\limits_{n=0}^{\infty}c_ns^n,\\
G_0(s)&=1+\left(1-\frac{1}{\varphi}\right)s-\frac{s^2}{\varphi}+2(1-\varphi)s^3-s^4+\frac{-9+5\sqrt{5}}{2}s^6+\ldots=\sum_{n=0}^{\infty}d_ns^n,
\end{align*}
where
\begin{align*}
&c_0 = 1,\, c_1=1-\varphi,\, c_2 = -1, \, c_3 = -\frac{1}{\varphi},\, c_n = \frac{\varphi^2 c_{n-1}-\varphi c_{n-2}-c_{n-3}}{\varphi^2},\, n = 4,\, 5,\, \ldots, \\
&d_0=1,\,d_1=1-\frac{1}{\varphi},\,d_2=-\frac{1}{\varphi},\,d_n=\frac{\varphi^2d_{n-1}-\varphi d_{n-2}-d_{n-3}}{\varphi^2},\,n=3,\,4,\,\ldots
\end{align*}
and the radius of convergence of both series $M_0(s)$ and $G_0(s)$ is $R_M=R_G\approx0.95747$; this is the modulus of $\approx0.61888 \pm 0.73057 i$, which is the roots of the smallest modulus of $1/(1-s)/G\alpha_0(s)$. Thus, neither series converges as $s\to1$.

According to eq. \eqref{alpha_ne_diff} and eq. \eqref{R_1}, 
\begin{align*}
{\pmb \alpha}_0&=\left\{1,\,\frac{3-\sqrt{5}}{2},\,\frac{1-\sqrt{5}}{2},\,1-\sqrt{5},\,-1,\,0,\,\frac{-9+5\sqrt{5}}{2},\,\ldots\right\}, \\
&\approx\{1,\,0.38197,\,-0.61803,\,-1.23607,\,-1,\,0,\,1.09017,\,\ldots\},\\
&\limsup_{n\to\infty}|\alpha_0(n)|^{1/n}=1/R_M\approx1.04442.
\end{align*}  

\subsection{Some famous sequences}

Propositions \ref{ex_zeta_1}--\ref{proposition_pi} show that the sequence ${\pmb b}$ can be chosen in a way that ${\pmb a}$, computed by \eqref{eq:seq_v1}, or $\alpha$-sequences, computed by \eqref{rec_for_expres}, give a desired sequence. In this section, we provide more sequences ${\pmb b}$ that determine some famous sequences ${\pmb a}={\pmb \alpha}_0$, i.e. we set $m=1$ and $a_0=1$:
\begin{enumerate}
    \item If $b_0 = 1,\, b_1 = 1/2, \, b_n = -5/2^n,\, n=2,\,3,\,\ldots$, then
    $2{\pmb a} = \{2,\,1,\,3,\,4,\,7,\,11,\,18,\,\ldots\} $ are Lucas numbers \cite[\href{https://oeis.org/A000032}{A000032}]{oeis}.
    \item If $b_0 = 1,\, b_1 = -6,\, b_2 = 26,\, b_3 = -84,\, b_n = -2b_{n-1}+4b_{n-2},\, n=4,\,5,\,\ldots$, then ${\pmb a} = \{1,\,7,\,23,\,63,\,159,\,383,\,895,\,\ldots\}$ are Woodall's (or Riesel's) numbers \cite[\href{https://oeis.org/A003261}{A003261}]{oeis}.
    \item If $b_0= b_3 =1,\,b_1=-1,\,b_2=0,\,b_n = -b_{n-3},\,n=4,\,5,\,\ldots$, then \\ ${\pmb a} = \{1,\,2,\,4,\,7,\,11,\,16,\,22,\,\ldots\} $ is Lazy Caterer's sequence \cite[\href{https://oeis.org/A000124}{A000124}]{oeis}.
    \item If ${\pmb b} = \{1,\,-1,\,-1,\,0,\,0,\,\ldots\}$, then ${\pmb a} = \{1,\,2,\,5,\,12,\,29,\,70,\,169,\,\ldots\}$ are shifted Pell's numbers \cite[\href{https://oeis.org/A000129}{A000129}]{oeis}.
    \item If ${\pmb b} = \{1,\,-1,\,1,\,0,\,0,\,\ldots\}$, then ${\pmb a} = \{1,\,2,\,3,\,4,\,5,\,6,\,7,\,\ldots\}$ are natural numbers \cite[\href{https://oeis.org/A000027}{A000027}]{oeis}.
    \item If ${\pmb b} = \{1,\,0,\,-1,\,-2,\,-6,\,-22,\,\ldots\}$, where $b_n,\, n\geqslant2$ are negative numbers of a series \cite[\href{https://oeis.org/A074664}{A074664}]{oeis} excluding its first value, then ${\pmb a} = \{1,\,1,\,2,\,5,\,15,\,52,\,203,\,\ldots\}$ are  Bell's numbers \cite[\href{https://oeis.org/A000110}{A000110}]{oeis}.
    \item If ${\pmb b} = \{1,\,0,\,-1,\,-2,\,-5,\,-14,\,\ldots\}$, where $b_n,\, n\geqslant2$ are negative Catalan's numbers excluding its first value, then ${\pmb a} = \{1,\,1,\,2,\,5,\,14,\,42,\,132,\,\ldots\}$ are  Catalan's numbers \cite[\href{https://oeis.org/A000108}{A000108}]{oeis}.
    \item If ${\pmb b} = \{1,\,1,\,-1,\,-2,\,-5,\,-14,\,\ldots\}$, where $b_n,\, n\geqslant2$ are negative Catalan's numbers excluding its first value, then ${\pmb a} = \{1,\,0,\,1,\,2,\,6,\,18,\,57,\,\ldots\}$ are shifted Fine's numbers \cite[\href{https://oeis.org/A000957}{A000957}]{oeis}.
    \item If $b_0 = 1,\, b_1 = b_3 = -1,\, b_2 = 0,\, b_4 = 3,\, b_n = b_{n-2}-b_{n-3},\, n=5,\,6,\,\ldots$, then ${\pmb a} = \{1,\,2,\,4,\,9,\,17,\,33,\,61,\,\ldots\}$ is Les Marvin's sequence \cite[\href{https://oeis.org/A007502}{A007502}]{oeis}.
    \item If ${\pmb b} = \{1,\,-1,\,2,\,-4,\,10,\,-20,\,\ldots\}$, where $b_n,\, n\geqslant2$ are numbers of a series \cite[\href{https://oeis.org/A308986}{A308986}]{oeis} excluding its first two values then ${\pmb a} = \{1,\,2,\,2,\,4,\,2,\,4,\,4,\,\ldots\}$ is Gould's sequence \cite[\href{https://oeis.org/A001316}{A001316}]{oeis}.
    \item If ${\pmb b} = \{1,\,0,\,-1,\,-1,\,-2,\,-4,\,-9,\,\ldots\}$, where $b_n,\, n\geqslant2$ are negative Motzkin's numbers, then ${\pmb a} = \{1,\,1,\,2,\,4,\,9,\,21,\,51,\,\ldots\}$ are Motzkin's numbers \cite[\href{https://oeis.org/A001006}{A001006}]{oeis}.
    \item If $b_0 = b_1 = 1,\, b_{2n} = 0,\, b_{2n+1} = -1,\,n\in\mathbb{N}$, then ${\pmb a} = \{1,\,0,\,0,\,1,\,0,\,1,\,1,\,\ldots\}$ are Padovan's numbers \cite[\href{https://oeis.org/A000931}{A000931}]{oeis}.
    \item If ${\pmb b} = \{1,\,25,\,324,\,3200,\,25650,\,\ldots\}$, where $b_n,\, n\geqslant2$ are numbers from a series \cite[\href{https://oeis.org/A006922}{A006922}]{oeis} excluding its first two values, then\\ ${\pmb a} = \{1,\,-24,\,252,\,-1472,\,4830,\,-6048,\,-16744,\,\ldots\}$ are Ramanujan's numbers \cite[\href{https://oeis.org/A000594}{A000594}]{oeis}.
    \item If ${\pmb b}=\{1,\,1-k,\,0,\,0,\,0,\,0,\,\ldots\}$, where $k\in\mathbb{N}_0$, then ${\pmb a}=\{k^n,\,n\in\mathbb{N}_0\}$.
    \item If ${\pmb b}=\{1,\,0,\,\underbrace{-1,\,-1,\,\ldots,\,-1}_{k},\,0,\,0,\,\ldots\}$, where $k$ is a number of "$-1$", then ${\pmb a}$ is $k$-generalized Fibonacci sequence \cite{Noe2005}.
 \end{enumerate}

Of course, this section of examples, including Propositions \ref{ex_zeta_1}--\ref{proposition_pi}, can be further extended, and some examples can (must?) be studied more deeply.

\section{Acknowledgments}

We thank James Tuite for pointing out that the root of the largest modulus of $1/(2s^3)+s/2=1$ from \cite[Example 1]{AA} coincides with the limit $T(n+1)/T(n),\,n\to\infty$, where $T(n)$ is the Tribonacci sequence. This inspired the presence and usage of Corollary \ref{Tuite}.

\bibliographystyle{plain}
\bibliography{bibliography}

\end{document}